\pdfoutput=1
\documentclass[final,leqno,onefignum,onetabnum,onealgnum,oneeqnum,onethmnum]{siamart171218}
\usepackage{lipsum}
\usepackage{amsfonts}
\usepackage{mathrsfs}
\usepackage{amsmath}
\usepackage{amssymb}
\usepackage{graphicx}
\usepackage{epstopdf}
\usepackage[caption=false]{subfig}
\graphicspath{{imgs/}}
\usepackage[algo2e,ruled,vlined,resetcount,linesnumbered,algosection]{algorithm2e}
\crefname{algocf}{alg.}{algs.}
\Crefname{algocf}{Algorithm}{Algorithms}
\Crefname{Algorithm}{Algorithm}{Algorithms}
\crefformat{figure}{#2{}\scshape{Fig.} #1{}#3}
\crefformat{lemma}{#2{}\scshape{Lemma} #1{}#3}
\crefformat{theorem}{#2{}\scshape{Theorem} #1{}#3}
\crefformat{proposition}{#2{}\scshape{Proposition} #1{}#3}
\crefformat{definition}{#2{}\scshape{Definition} #1{}#3}
\newtheorem{remark}{Remark}
\crefformat{remark}{#2{}\scshape{Remark} #1{}#3}
\AppendGraphicsExtensions{.PNG}
\usepackage{etoolbox}
\apptocmd{\sloppy}{\hbadness 10000\relax}{}{}
\makeatletter
\newcommand{\clearsubcaptcounter}{\setcounter{sub\@captype}{0}}
\makeatother
\usepackage{siunitx}
\def\Lbd{\boldsymbol{\Lambda}}
\def\lbd{\lambda}
\def\mbf#1{\mathbf{#1}}
\def\D{\nabla}
\def\a{\alpha}
\def\te{\theta}
\def\ha{\hat{\a}}
\def\hth{\hat{\te}}
\def\O{\Omega}
\def\G{\Gamma}
\def\S{\Sigma}
\def\s{\sigma}
\usepackage{tikz}
\usetikzlibrary{calc}

\title{Sparse-data based 3D surface reconstruction with vector matching}
\author{
    Bin Wu\thanks{Department of Computer Science, Electrical Engineering and Mathematical Sciences, Western Norway University of Applied Sciences, Inndalsveien 28, 5063 Bergen, Norway (\email{bin.wu@hvl.no}, \email{talal.rahman@hvl.no}).}
    \and
    Xue-Cheng Tai\thanks{Department of Mathematics, University of Bergen, All\'{e}gaten 41, 5007 Bergen, Norway (\email{xue-cheng.tai@uib.no}).}
    \and
    Talal Rahman\footnotemark[1]
}
\headers{Surface reconstruction with vector matching}{B. Wu, X. C. Tai, and T. Rahman}
\begin{document}
\maketitle

\begin{abstract}
    Three dimensional surface reconstruction based on two dimensional sparse information in the form of only a small number of level lines of the surface with moderately complex structures, containing both structured and unstructured geometries, is considered in this paper. A new model has been proposed which is based on the idea of using normal vector matching combined with a first order and a second order total variation regularizers. A fast algorithm based on the augmented Lagrangian is also proposed. Numerical experiments are provided showing the effectiveness of the model and the algorithm in reconstructing surfaces with  detailed features and  complex structures for both synthetic and real world digital maps.
\end{abstract}

\begin{keywords} 
    total variation regularization, surface reconstruction, augmented Lagrangian
\end{keywords}

\begin{AMS}
    65F05
\end{AMS}

\section{Introduction}

Reconstructing an image from its little information available is a challenging but interesting task for image processing, which has attracted much attention over the years, finding applications in many areas where, for instance, information are available only in the form of a small number of level lines or isolated points. Two of its well known applications are the sketch based design where three dimensional figures or structures are designed from sketches made by an artist or a computer program, and the surface reconstruction where three dimensional surface is reconstructed from its level lines, see for instance two recent publications \cite{HQSJTS2013, LMS2013}, for a review on the subject. In this paper, we are interested surfaces that are moderately complex, in the sense that they contain both structured and unstructured geometries, man made or natural, and feature both kinks (sharp twists) and creases (folds). We are particularly interested in a variational model that combines both the features of a sketch based design and that of a surface reconstruction, by incorporating both the height and the vector information from the level lines in a natural way so that a more precise design or reconstruction is possible. 

Sketch based design has been a popular way of design in three dimensions, cf. e.g. \cite{ZHH1996, IMH1999, KH2006, NISA2007, OSSJ2009}, because it is both intuitive and effective, particularly in the games and cartoon design. In a sketch based design, the information is available either in the form of contour lines (level lines) with or without the height values, cf. e.g. \cite{IMH1999}, complex sketches with elevations, cf. e.g. \cite{KH2006}, or structured annotations, cf. e.g. \cite{GIZ2009}. The class of algorithms presented in these papers are however limited in their capabilities, particularly when it comes to reconstructing structures with crease. Although there exists a work on artificially adding crease to the design, cf. e.g. \cite{OSSJ2005}, for large and complex images, such algorithms become less effective and computationally more expensive. A different class of models, based on the total variation minimization, turned out to be much more effective, cf. \cite{HQSJTS2013}. The model presented in these two papers is based on interpolating normal vectors under the curl-free constraint and then reconstructing the 3D surface from the obtained vector field. Inspired by the use of surface gradients in surface reconstruction, cf. e.g. \cite{ARC2006, FCM1987, PCFKH2001, SCSJ1990, ZDPS2002, WTBS2007, PF2006, NWT2009}, the model is an extension of the TV-Stokes model, cf. \cite{RTO2007}, to surface reconstruction. The TV-Stokes models are based on using the curl-free or the divergence free constraint in image processing, cf. e.g. \cite{EMR2009, TBH2009, HTBB2011, LRT2011} for further works on TV-Stokes models. The model of \cite{HQSJTS2013} performs very well in preserving both edges and crease structures, however, it requires prior information on the vectors, e.g. the length of vectors \cite{WRT2017, HQSJTS2013}, which are not always available.

Three dimensional (3D) surface reconstruction from contours or isolated points with height values, has been the second most popular application area. Unlike the three dimensional design, the height values here are needed because the reconstructed surfaces are expected to be as close as possible to the ground truth, e.g., in digital elevation maps and data compression, cf. e.g. \cite{LMS2013}. One approach to solve the problem is to use explicit parametrization of the given contours, with subsequent point wise matching and interpolation between the contours, cf. e.g. \cite{MSS1992, MM1998, M2011}. In some cases, an explicit parametrization may be difficult and expensive to compute, and a loss of continuity of slope across contours may become a challenge to deal with. An alternative way is to consider the surface as a function over the domain, and interpolate the function based on solving partial differential equations, as in the AMLE (Absolutely Minimizing Lipschitz Extension) model, cf. e.g. \cite{AGLM1993, CMS1998}. Although AMLE interpolation is able to interpolate data given on isolated points and on level lines, it has the drawback that the level lines of its interpolants are smooth making it difficult to preserve kinks, and consequently creases on the surface, plus it cannot interpolate slopes of the surface. To overcome this, one has to rely on higher-order methods or regularizers, cf. e.g. \cite{F1982, M1984, CFB1997, MM1999, LMS2013}. 

Particularly interesting to our present work, is the model proposed in \cite{LMS2013}, a variational model, which uses a third order anisotropic regularizer whose anisotropy is based on an auxiliary vector field connecting adjacent level lines; it is the field of directions formed explicitly in a separate step, in which the normals of the level lines change the least. The model performs very well for surface reconstruction, particularly in preserving the geometry of the given level lines, and propagating it smoothly across the level lines. The model does not require any regularity of the level lines, however, it requires the level lines lying next to each other to be similar, so that points lying on them can be associated. This may not always be the case, in particular when there are only few level lines available.

Our aim is to recover the surface from the few level lines that are available, and somehow use the geometry of those level lines in a physically consistent way providing further precision to our reconstruction. In case of a 3D design, which we will not cover in this paper, these level lines will be the sketch lines drawn by an artist. We propose a simple one step variational model incorporating both the height and the vector information in the model, the vectors being the unit normal vectors, either calculated from the level lines itself or provided by the artist on the sketches. The model consists of a first-order and a second-order total variation (isotropic) regularizer under a fidelity constraint on the height, together with a vector (normal vector) matching term to account for the anisotropy in the model.

We do not impose any regularity on the level lines, nor do we assume any similarity between the level lines close to each other. In case of very sparse data this is very likely to happen. We require the geometry of the level lines, in particular, the non-differentiability along level lines, as well as the smoothness of the gradient across level lines be preserved in the interpolated surface. The model allows for adaptive adjustment of the normal vectors, consequently the shape of the surface, and an almost perfect reconstruction even with a small number of level lines for complex surface with mixed geometries. We propose a fast algorithm for the numerical solution, which is based on the augmented Lagrangian method \cite{GLT1989, TW2009, WZT2011}, featuring sub-problems with closed form solutions and fast iterative solution.

The paper is organized as follows. In section 2, we propose our model for surface reconstruction from a set of sketches or level lines, and individual points, in section 3, we present a fast algorithm for the numerical solution, based on the augmented Lagrangian, and in section 4, we present our numerical experiments on both synthetic data and real data for the verification. Finally, in section 5, we give our conclusion.

\section{The proposed model}
\label{sec:proposing}

\begin{figure}[!htbp]
    \captionsetup[subfigure]{
                            justification=centering,
                            }
    \centering
    \begin{minipage}{1.0\textwidth}
        \centering
        \subfloat
            {
            \includegraphics[width=0.6\linewidth,
                            angle=0]{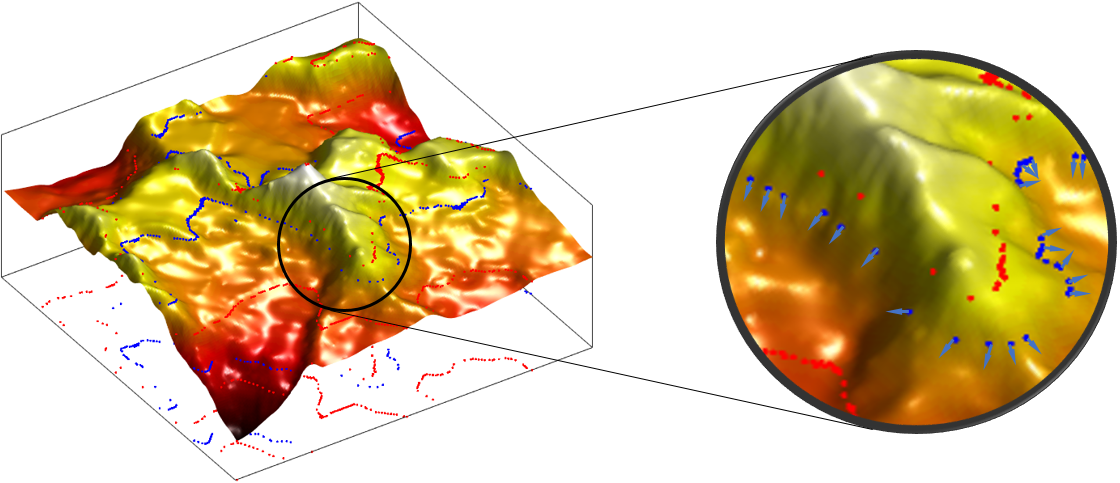}
            \label{subfig:dm}
            }
        \vspace{-0.01\linewidth}
    \end{minipage}
    \caption{Given 3 dimensional surface on a 2 dimensional domain. $\G$ (the set of blue points) corresponds to the points where unit vectors are given; $\S$ (the set of red points) corresponds to the points where the elevation/height values are given. These are points that are either isolated or given on level lines.}
    \label{fig:domain}
\end{figure}

The problem we consider is an inverse problem to recover the two dimensional (2D) height map $I:\O \mapsto \mathbb{R}$ on a domain $\O \subset \mathbb{R}^2$ in two space dimensions, from known sparse data given in the form of level lines as $C_i = \{ x: I(x) = I_{\S i} \}$ corresponding to the given elevation $I_{\S i}$, for $i=1,\ldots,N$. The collection of all the given level lines is denoted by $\S$, that is $C_i \in \S$. The aim is to propose a model that is based on the total variation minimization and an explicit use of the geometry of the level lines in the functional. The use of this geometry, in one form or another, is essential for the surface reconstruction, e.g., for capturing kinks, specially when the level lines have coastal like structures. However, how well the kinks are represented depends on how well the geometry of the level lines has been incorporated into the model. There are models which include such information, cf. e.g. \cite{HQSJTS2013, LMS2013}. However, the use of such information, in these papers, has been done somewhat indirectly, that is by first constructing a smooth vector field from the vector information available on the level lines, that is in a prepossessing step, and then using the new vector field in the reconstruction step. These models are quite powerful and effective, however, there are cases where these models fall short; we refer to the numerical section of this paper for one such case. We propose a more direct approach, where we introduce a vector matching term into our model, involving explicitly the vector information available on the level lines, as follows:
\begin{equation}
    \min_I  
    \left \{
        g | \D (\D I) |_{F \ (\O)} 
        + h | \D I |_{ \ (\O)} 
        + \a V(\D I, \mbf{v}_{\G})_{ \ (\G)} 
        + {\te} | I - I_{\S} |_{2 \ (\S)}^2
    \right \},
    \label{eqs:opt}
\end{equation} 
where the first two terms are the second order and the first order total variational regularizer, respectively, the third term is the vector matching term, and the fourth term is the data fidelity term. $| \cdot |_{F \ (\O)}$, $| \cdot |_{ \ (\O)}$, and $| \cdot |_{2 \ (\S)}$ are the Frobenius-, $L^1$-, and $L^2$-, norm, respectively, over the domain specified inside bracket. The vector matching term is the integral over $\G$ of function $V$. $g$, $h$, $\a$ and $\te$ are the scalar parameters. $\G$ is the set of points where the unit vectors $\mbf{v}_{\G}$ are given and $\S$ is the set of points where the elevation/height values $I_{\S}$ are given, see \cref{fig:domain} for an illustration. We note that the vector $\mbf{v}_{\G}$ is either given, or can be extracted from the given level lines.

\subsubsection*{The vector matching term}
    On a given level line the vector $\mbf{v}_{\G}$ can either be given as the unit tangent vector or be given as the unit normal vector. The idea is to match $\D I^*$, that is the gradient of the reconstructed surface $I^*$, to this vector in some consistent way, e.g., $\D I^*$ should be orthogonal to $\mbf{v}_{\G}$ if $\mbf{v}_{\G}$ is the unit tangent vector and parallel otherwise.

    In case of the unit tangent vector, the only way this matching can happen is through minimizing a norm of $\D I \cdot \mbf{v}_{\G}$, e.g.,
    \begin{equation}
        V(\D I,\mbf{v}_{\G})_{ \ (\G)} 
        := 
        | \D I \cdot \mbf{v}_{\G} |_{2 \ (\G)}^2.
        \label{eqs:tangent}
    \end{equation}
    If we use this norm, the corresponding minimization problem becomes as follows,
    \begin{equation}
        \min_I  
        \left \{ 
            g | \D (\D I) |_{F \ (\O)} 
            + h | \D I |_{ \ (\O)} 
            + \a | \D I \cdot \mbf{v}_{\G} |_{2 \ (\G)}^2 
            + {\te} | I - I_{\S} |_{2 \ (\S)}^2
        \right \}.
        \label{eqs:opt_tangent}
    \end{equation}
    Its minimum is attained when $\D I^*$ ($\D I^* \neq \mbf{0}$) is perpendicular to the given vector $\mbf{v}_{\G}$, or $\D I^* = \mbf{0}$. There are basically two shapes that the minimization generates. In the first case, when $\D I^* \perp \mbf{v}_{\G}$, $\D I^* \neq \mbf{0}$, $\a$ can take any number, also when $\a=0$, which is an isotropic case. In the second case, when $\D I^* = \mbf{0}$, the minimization generates flattened structures.

    In case of the unit normal vector, the only consistent way is minimizing the inner product $ - \D I \cdot \mbf{v}_{\G \ (\G)}$, that is
    \begin{equation}
        V(\D I,\mbf{v}_{\G})_{ \ (\G)} 
        := 
        - \D I \cdot \mbf{v}_{\G \ (\G)}.
        \label{eqs:normal}
    \end{equation}
    With this, the corresponding minimization problem becomes as follows,
    \begin{equation}
        \min_I  
        \left \{ 
            g | \D (\D I) |_{F \ (\O)} 
            + h | \D I |_{ \ (\O)} 
            - \a \D I \cdot \mbf{v}_{\G \ (\G)} 
            + {\te} | I - I_{\S} |_{2 \ (\S)}^2
        \right \}.
        \label{eqs:opt_normal}
    \end{equation}
    Its minimum is attained only when $\D I^*$ is pointing in the same directions as the unit vector $\mbf{v}_{\G}$.
    With this matching term in the functional, cf. (\ref{eqs:opt_normal}), any irregularities in the given level lines will be preserved, which is then stretched smoothly across the level lines by the regularizers. Also, by varying the parameter $\a$ we can control the shape of the reconstructed surface, see the numerical section for an illustration. For our model, we therefore choose $\mbf{v}_{\G}$ to be the unit normal vector on the level lines.

\subsubsection*{The regularizer terms}
    To explain our choice of regularizer, we consider the following single regularizer model, where we consider only one regularizer term in the functional, with varying order, that is
    \begin{equation}
        \min_I 
        \left \{
            g | \D^k I |_{F \ (\O)} 
            - \a \D I \cdot \mbf{v}_{\G \ (\G)} 
            + \te | I - I_{\S} |_{2 \ (\S)}^2
        \right \},
        \label{eqs:opt_k_ord}
    \end{equation}
    for $k=1,2,\ldots$. The notations are the same as before, cf. (\ref{eqs:opt_normal}). Choosing $\a$ small compared to both $g$ and $\te$, the optimal solution of the problem (\ref{eqs:opt_k_ord}) at any position away from $\S$ and $\G$ would give $\D ^k I \sim 0$. For $k=1$, it corresponds to $\D ^1 I \sim 0$, in other words no change in the height value, thereby producing a staircase effect on the surface. For $k=2$, this corresponds to the Hessian $\D ^2 I \sim 0$, implying no change in the gradient of the height, thereby preserving the slope, consequently the shape, of the surface. For our model, it is enough to use the two, or a combination of the two. Increasing the number $k$ introduces oscillations to the solution, becoming unbounded near the boundary, see the numerical section for an illustration.  

    Our model is a combination of a first order and a second order isotropic regularizer with an anisotropic matching term, which is different from the anisotropic third order regularizer of \cite{LMS2013}, in the sense that the anisotropy in our model is separated from the regularizer, while in the latter the anisotropy is built in the regularizer itself.

\subsubsection*{The choice of parameters} The choice of parameters $g$, $h$, $\a$ and $\te$ is crucial for our model In general the data fidelity parameter $\te$ is set to be large enough to ensure the reconstructed surface is consistent with the given data. The choice of vector matching parameter $\a$ depends on the surface to be reconstructed, for instance, the bigger the value the steeper the structure becomes. The parameter $h$ (for the first order total variation) is used for flat or almost flat structure. The larger its value the stronger the flattening is. The parameter $g$ (for the second order total variation) is used for sloped structure smooothly connecting the level lines across. The larger the $g$ the smoother the surface between the level lines becomes.
    
\section{The numerical algorithm}
\label{sec:alg}

In this section, we introduce our algorithm based on the augmented Lagrangian for solving the optimization problem (\ref{eqs:opt}), cf. \cite{GLT1989, NW2006} for general literature on augmented Lagrangian method, and \cite{TW2009, WZT2011} for its applications in image processing.

We denote $\ha$ and $\hth$ as follows,
\begin{equation*}
    \ha 
    =
    \begin{cases}
        \a,             \hspace{1em}    
        \mbox{on}       \hspace{0.5em}  
        \G              \\
        0,              \hspace{1em}    
        \mbox{in}       \hspace{0.5em}  
        \O \backslash \G,
    \end{cases} 
    \quad \mbox{and} \quad
    \hth 
    =
    \begin{cases}
        \te,            \hspace{1em}    
        \mbox{on}       \hspace{0.5em}  
        \S              \\
        0,              \hspace{1em}    
        \mbox{in}       \hspace{0.5em}  
        \O \backslash \S.
    \end{cases} 
\end{equation*}
The problem (\ref{eqs:opt}) is thus defined on the whole of domain $\O$ as, 
\begin{equation}\label{eqs:opt_hat}
    \min_I  
    \left \{ 
        g | \D (\D I) |_{F} 
        + h | \D I | 
        - \ha \D I \cdot \mbf{v}_{\G} 
        + \hth | I - I_{\S} |_{2}^2 
    \right \}.  
\end{equation}

Accordingly, we introduce auxiliary variables for the derivatives in $L^{1}$-terms turning the unconstrained optimization problem into a constrained optimization problem, cf. e.g., \cite{NW2006,TW2009,WZT2011}.

We set $\mbf{Q} := \D \mbf{E}$, $\mbf{P} := \D I$ and $\mbf{E} := \mbf{P}$, where $\mbf{E} \in \mathbb{R}^2$ and $\mbf{P} \in \mathbb{R}^2$ are 2-dimensional vectors, and $\mbf{Q} \in \mathbb{R}^{2\times2}$ is a 2-by-2 matrix. Using the new variables, we get the following constrained minimization problem.
\begin{equation*}
    \min_{\mbf{Q},\mbf{P},I} 
    \left \{ 
        g | \mbf{Q} |_{F} 
        + h | \mbf{P} | 
        - \ha \mbf{P} \cdot \mbf{v}_{\G} 
        + \hth | I - I_{\S} |_{2}^2
    \right \},
    \end{equation*} 
subject to
\begin{equation*}
    \mbf{P} - \D I = 0,      ~\ 
    \mbf{E} - \mbf{P} = 0,    ~\ 
    \mbox{and}                      ~\ 
    \mbf{Q} - \D \mbf{E} = 0.
\end{equation*} 
Assigning to each constraint a Lagrange multiplier and a penalty term, the Lagrangian functional reads as follows,
\begin{eqnarray}
    \mathscr{L}
        (
            \mbf{Q}, \mbf{P}, \mbf{E}, I; 
            \Lbd_Q,  \Lbd_P,  \Lbd_E
        )~
    &=&  
        ~g | \mbf{Q} |_F 
        ~+~ h | \mbf{P} | 
        ~-~ \ha \mbf{P} \cdot \mbf{v}_{\G}
        ~+~ \hth | I - I_{\S} |^2 
        \\ 
        && 
        +~ \Lbd_Q \cdot (\mbf{Q} - \D \mbf{E}) 
        ~+~ \cfrac{c_Q}{2} | \mbf{Q} - \D \mbf{E} |_F ^2  
        \nonumber \\ 
        && 
        +~ \Lbd_P \cdot (\mbf{P} - \D I) 
        ~+~ \cfrac{c_P}{2} | \mbf{P} - \D I | ^2
        \nonumber \\
        && 
        +~ \Lbd_E \cdot (\mbf{E} - \mbf{P}) 
        ~+~ \cfrac{c_E}{2} | \mbf{E} - \mbf{P} | ^2,
        \nonumber 
\end{eqnarray}
where $\Lbd_Q$, $\Lbd_P$, and $\Lbd_E$ are the Lagrange multipliers, $c_Q$, $c_P$, and $c_E$ are the positive penalty parameters. The augmented Lagrangian method is to seek the saddle point of the following problem:
\begin{equation}\label{eqs:opt_saddle}
    \min_{\mbf{Q}, \mbf{P}, \mbf{E}, I}    \   \
    \max_{\Lbd_Q, \Lbd_P, \Lbd_E} 
        \mathscr{L}
            (
                \mbf{Q}, \mbf{P}, \mbf{E}, I;
                \Lbd_Q,  \Lbd_P,  \Lbd_E
            ).
\end{equation}
For the solution, we solve its associated system of optimality conditions with an iterative procedure, see \Cref{alm:algorithm1} and \Cref{alm:algorithm2}. For the convenience, we use $\Lbd := (\Lbd_{Q}, \Lbd_{P}, \Lbd_E)$ to denote the set of Lagrange multipliers.

\begin{algorithm2e}[!htbp] 
    \linespread{1.2}
    \selectfont
    \SetAlgoLined
    \KwIn{$\mbf{v}_{\G}$ (The unit vector along $\G$) and $I_{\S}$ (The elevation along $\S$)} 
    \KwOut{$I$ (The elevation on the whole domain $\O$)}
    Set $k = 0$; \\
    Initialize each of $\mbf{Q}^0$, $\mbf{P}^0$, $\mbf{E}^0$, $I^0$ and $\Lbd^0$ to be zero\;
    \While{not converged}
        {
            Update k = k + 1 \; 
            Solve for $\mbf{Q}, \mbf{P}, \mbf{E}$ and $I$: 
            \vspace{-1em}
            \begin{equation}\label{eqs:opt_min}
                (\mbf{Q}^k, \mbf{P}^k, \mbf{E}^k, I^k ) 
                    = \mbox{arg} \min \limits_{\mbf{Q},\mbf{P},\mbf{E},I} \mathscr{L}  (\mbf{Q},\mbf{P},\mbf{E},I;\Lbd^{k-1}) ;
            \end{equation} \\
            \vspace{-1.5em}
            Update $\Lbd = (\Lbd_Q, \Lbd_P, \Lbd_E)$:\\			
                    $ \qquad \qquad \qquad
                      \Lbd^k_Q = \Lbd^{k-1}_Q + c_Q(\mbf{Q}^k - \D \mbf{E}^k), 
                    $ \\
                    $ \qquad \qquad \qquad
                      \Lbd^k_P = \Lbd^{k-1}_P + c_P(\mbf{P}^k - \D I^k), 
                    $ \\
                    $ \qquad \qquad \qquad
                      \Lbd^k_E = \Lbd^{k-1}_E + c_E(\mbf{E}^k -\mbf{P}^k); 
                    $
        }
    Return
        $I = I^k$;
    \caption{The augmented Lagrangian for the problem (\ref{eqs:opt_saddle})} 
    \label{alm:algorithm1}       
\end{algorithm2e}         

\vspace{.2cm}
\noindent Because the variables $\mbf{Q}$, $\mbf{P}$, $\mbf{E}$ and $I$ in $\mathscr{L} (\mbf{Q},\mbf{P},\mbf{E},I;\Lbd^{k-1})$ are coupled together in the minimization problem (\ref{eqs:opt_min}), it is difficult to solve them simultaneously. We split the minimization problem into four sub minimization problems, and solve them alternatingly until convergence, cf. \Cref{alm:algorithm2}. Typically we need only one iteration ($L=1$). 

\begin{algorithm2e}[!htbp]
    \linespread{1.2}
    \selectfont
    \SetAlgoLined
    \KwIn{$\mbf{Q}^{k-1}$, $\mbf{P}^{k-1}$, $\mbf{E}^{k-1}$, $I^{k-1}$ and $\Lbd^{k-1}$}
    \KwOut{$\mbf{Q}^k$, $\mbf{P}^k$, $\mbf{E}^k$ and $I^k$}
    Set $l = 0$; \\
    Initialize
    $\mbf{Q}^{k,0} = \mbf{Q}^{k-1}$, $ \mbf{P}^{k,0} = \mbf{P}^{k-1}$, $ \mbf{E}^{k,0} = \mbf{E}^{k-1}$ and $I^{k,0} = I^{k-1}$; \\
    \While{not converged and $l < L$ }
        {   Solve for $\mbf{Q}$ (the $\mbf{Q}$-subproblem) \\
            $ \qquad \qquad \qquad
                \mbf{Q}^{k,l +1} 
                = \mbox{arg} \min \limits_{\mbf{Q}} 
                    \mathscr{L}(
                        \mbf{Q},\mbf{P}^{k,l},\mbf{E}^{k,l},I^{k,l}; \Lbd^{k-1}
                                );
            $ \\
            Solve for $\mbf{P}$ (the $\mbf{P}$-subproblem) \\
            $ \qquad \qquad \qquad
                \mbf{P}^{k,l +1} 
                = \mbox{arg} \min \limits_{\mbf{P}}
                    \mathscr{L}(
                        \mbf{Q}^{k,l+1},\mbf{P},\mbf{E}^{k,l},I^{k,l}; \Lbd^{k-1}
                                );
            $ \\
            Solve for $\mbf{E}$ (the $\mbf{E}$-subproblem) \\
            $ \qquad \qquad \qquad
                \mbf{E}^{k,l +1} 
                = \mbox{arg} \min \limits_{\mbf{E}}
                    \mathscr{L}(
                    \mbf{Q}^{k,l+1},\mbf{P}^{k,l+1},\mbf{E},I^{k,l}; \Lbd^{k-1}
                                );
            $ \\
            Solve for $I$ (the $I$-subproblem) \\
            $ \qquad \qquad \qquad
                I^{k,l +1} 
                = \mbox{arg} \min \limits_{I} 
                    \mathscr{L}( 
                        \mbf{Q}^{k,l+1},\mbf{P}^{k,l+1},\mbf{E}^{k,l+1},I; \Lbd^{k-1}
                                );
            $ \\
            Update $l = l + 1$\;
        }
    Return
    $\mbf{Q}^{k} = \mbf{Q}^{k,l}$, $ \mbf{P}^{k} = \mbf{P}^{k,l}$, $ \mbf{E}^{k} = \mbf{E}^{k,l}$ and $I^{k} = I^{k,l}$;
    \caption{Solve the minimization problem (\ref{eqs:opt_min}). }
    \label{alm:algorithm2}
\end{algorithm2e}

We state an important relation, cf. \cref{rem:compatible}, which is used in our algorithm to give us a simpler approach to solve the first two sub-minimization problems of \Cref{alm:algorithm2}. 

\begin{remark} \label{rem:compatible}
    If $A$ and $B$ are two matrices such that $A=\lbd B$ for some non-negative scalar $\lbd$, then we say that $A$ is compatible with $B$. It is easy to see that $A/|A|_F = B/|B|_F$.
\end{remark}

\vspace{1.0em}
\noindent \textbf{The Q-subproblem:} The first problem is to solve for $\mbf{Q}$, freezing the other variables, in (\ref{eqs:opt_min}).
    \begin{equation}\label{eqs:sub_q}
        \mbf{Q}^{*} 
        = 
        \mbox{arg}\min_\mbf{Q} 
            \left \{ 
                g | \mbf{Q} |_F 
                + \Lbd_Q \cdot \mbf{Q}  
                + \cfrac{c_Q}{2} | \mbf{Q} - \D \mbf{E} |^2_F 
            \right \}. 
    \end{equation}
    We will find a closed form solution to this.
    The corresponding optimality condition, or the Euler-Lagrange equation, is as follows
    \begin{equation*}
        \cfrac{g}{c_Q} \cfrac{\mbf{Q}^{*} }{| \mbf{Q}^{*} |_F} 
        + \mbf{Q}^{*}  
        = 
        \D \mbf{E} - \cfrac{\Lbd_Q}{c_Q}.
    \end{equation*} 
    Since $g$ and $c_Q$ are both positive numbers and $| \mbf{Q}^{*} |_F$ is a positive number, the matrices $\mbf{Q}^{*}$ and $(\D \mbf{E} - {\Lbd_Q}/{c_Q})$ become compatible in the sense of \cref{rem:compatible} because 
    $
        ( {g}/{c_Q} \ |\mbf{Q}^{*}|_F^{-1} + 1 ) \mbf{Q}^{*}  
        = 
        \D \mbf{E} - {\Lbd_Q}/{c_Q}.
    $ According to which, we can replace the matrix 
    ${\mbf{Q}^{*}}/{| \mbf{Q} ^{*} |_F}$ with 
    $(\D \mbf{E} - {\Lbd_Q}/{c_Q})/| \D \mbf{E} - {\Lbd_Q}/{c_Q}|_F$ in the above equation. Moving the first term to the right hand side we get
    \begin{equation*}
        \mbf{Q}^{*} 
        = 
        \left(
                1 - \cfrac{g}{c_Q| \D \mbf{E} - \frac{\Lbd_Q}{c_Q} |_F}
        \right) 
        \left( 
            \D \mbf{E} - \cfrac{\Lbd_Q}{c_Q} 
        \right).
    \end{equation*} Again since we have already seen that in the sense of \cref{rem:compatible}, $\mbf{Q}^{*}$ and $(\D \mbf{E} - {\Lbd_Q}/{c_Q})$ are compatible, the coefficient $\left(1 - {g}/{c_Q} \ | \D \mbf{E} - {\Lbd_Q}/{c_Q}|_F^{-1} \right)$ must also be non-negative. Hence
    \begin{equation*}
        \mbf{Q}^{*}  
        = \max 
            \left\{
                0, 
                1 - \cfrac{g}{c_Q| \D \mbf{E} - \frac{\Lbd_Q}{c_Q} |_F} 
            \right\} 
        \left( \D \mbf{E} - \cfrac{\Lbd_Q}{c_Q} \right),
    \end{equation*}
    which is the solution to the $\mbf{Q}$-subproblem.

\vspace{1.0em}
\noindent \textbf{The P-subproblem:} The second problem is to solve for $\mbf{P}$, freezing the other variables, in (\ref{eqs:opt_min}).
    \begin{equation}\label{eqs:sub_p}
        \mbf{P}^{*} 
        = 
        \mbox{arg }\min_\mbf{P} 
            \left \{ 
                    h | \mbf{P} | 
                    - \ha \mbf{P} \cdot \mbf{v}_{\G} 
                    + (\Lbd_P - \Lbd_E) \cdot \mbf{P}  
                    + \cfrac{c_P}{2} | \mbf{P} - \D I |^2 
                    + \cfrac{c_E}{2} | \mbf{E} - \mbf{P} |^2 
            \right \}.
    \end{equation}
    We will find a closed form solution to this. The corresponding Euler-Lagrange equation is as follows
    \begin{equation*}
        \cfrac{h}{c_P + c_E} \cfrac{\mbf{P}^{*} }{| \mbf{P}^{*}  |} 
        + \mbf{P}^{*}  
        = 
        \mbf{X},
    \end{equation*} where $\mbf{X} := {(c_P \D I + c_E \mbf{E} - \Lbd_P + \Lbd_E + \ha \mbf{v}_{\G})}/{(c_P + c_E)}.$ Since $h$, $c_P$ and $c_E$ are positive numbers and $|\mbf{P}^{*}|$ is a positive number, the vectors $\mbf{P}^{*}$ and $\mbf{X}$ become compatible in the sense of \cref{rem:compatible} because 
    $
        ({h}/{(c_P + c_E)} \ |\mbf{P}^{*}|^{-1} + 1) \mbf{P}^{*} 
        = 
        \mbf{X}.
    $ According to which, we can replace the vector ${\mbf{P}^{*}}/{|\mbf{P}^{*}|}$ with $\mbf{X}/|\mbf{X}|$ in the above equation. Moving the first term to the right hand side we get
    \begin{equation*}
        \mbf{P}^{*}  
        = 
        \left( 
            1 - \cfrac{h}{(c_P + c_E)|\mbf{X}|} 
        \right)
        \mbf{X}.
    \end{equation*} Again since we have already seen that in the sense of \cref{rem:compatible}, $\mbf{P}^{*}$ and $\mbf{X}$ are compatible, the coefficient $(1 - {h}/{(c_P + c_E)} \ |\mbf{X}|^{-1})$ must also be non-negative. Hence
    \begin{equation*}
        \mbf{P}^{*}  
        = 
        \max
            \left\{
                0, 
                1 - \cfrac{h}{(c_P + c_E)| \mbf{X} |} 
            \right\} 
            \mbf{X},
    \end{equation*}
    which is the solution to to the $\mbf{P}$-subproblem.

\vspace{1.0em}
\noindent \textbf{The E-subproblem:} The third problem is to solve for $\mbf{E}$, freezing the other variables, in (\ref{eqs:opt_min}).
    \begin{equation}\label{eqs:sub_e}
            \mbf{E}^{*} 
            = 
            \mbox{arg}\min_\mbf{E} 
                \left \{ 
                        \Lbd_E \cdot \mbf{E} 
                        + \cfrac{c_E}{2} | \mbf{E} - \mbf{P} |^2 
                        - \Lbd_Q \cdot \D \mbf{E} 
                        + \cfrac{c_Q}{2} | \mbf{Q} - \D \mbf{E} |^2_F 
                \right \}.
    \end{equation}
    We will find a closed form solution to this. The corresponding Euler-Lagrange equation is the following.
    \begin{equation*}
        \D \cdot (\D \mbf{E}^{*})  
        - \frac{c_E}{c_Q} \mbf{E}^{*}  
        = 
        \frac{1}{c_Q} \Lbd_E 
        - \frac{c_E}{c_Q} \mbf{P} 
        + \frac{1}{c_Q}\D \cdot \Lbd_Q 
        + \D \cdot \mbf{Q},
    \end{equation*} which is a set of two inhomogeneous modified Helmholtz equations, one equation for each component of $\mbf{E} = (E_1,E_2)$, with the following Neumann boundary conditions, 
    \begin{equation*}
        \D E_1 \cdot \mbf{\nu} =  (\mbf{Q}_1 + \frac{1}{c_Q} \Lbd_{Q1}) \cdot \mbf{\nu},
         \quad \mbox{and} \quad
        \D E_2 \cdot \mbf{\nu} =  (\mbf{Q}_2 + \frac{1}{c_Q} \Lbd_{Q2}) \cdot \mbf{\nu},
    \end{equation*}
    where $\mbf{Q}_1$ and $\mbf{Q}_2$ are the row vectors of the matrix $\mbf{Q}$, and $\Lbd_{Q1}$ and $\Lbd_{Q2}$ are the corresponding Lagrange multipliers, respectively. $\mbf{\nu}$ is the outward unit normal vector on the boundary of the domain. Each equation is solved in the same way as follows.

    \vspace{1.0em}
    \noindent \textbf{Solving the inhomogeneous modified Helmholtz equation:}
    \begin{equation}\label{eqs:hlm}
        \triangle u(x,y) - \lbd u(x,y) = F(x,y),
    \end{equation}
    with a Neumann boundary condition and $\lbd$ a positive scalar, also known as the inhomogeneous modified Helmholtz equation. A fast solver based on discrete cosine transform similar for the Poisson equation, cf. \cite{VLC1992,EMR2009} and also cf. \cite{CAE2009} for details, is given below.

    The discrete version of Laplace operator in \cref{eqs:hlm} is the matrix
    \begin{equation*}
        \left[
            \begin{array}{cccccc}
             -1 &      1 &        &        &    \\
              1 &     -2 &  1     &        &    \\
                & \ddots & \ddots & \ddots &    \\
                &        &      1 &     -2 &  1 \\
                &        &        &      1 & -1 
            \end{array}
        \right],
    \end{equation*}
    with the help of discrete cosine transformation matrix $C \in \mathbb{R}^{N \times N}$ and singular value decomposition, we get the following decomposition,
    \begin{equation*}
        - C^{\top} 
        \left[ 
            \begin{array}{cc}
                0 &  \\
                  & \S ^2
            \end{array}
        \right] 
        C,
    \end{equation*}
    where $\S = \mbox{diag} (\s_1, \cdots, \s_{N−1})$ is the diagonal matrix with its entries representing the singular values $\s_k = 2 \sin ({\pi k}/{2}/N)$ for $k = 1,2,\cdots,N-1$. Substituting the decomposition back into \cref{eqs:hlm}, we get
    \begin{equation*}
        -u 
        C^{\top} 
        \left[ 
            \begin{array}{cc}
                0 &  \\
                  & \S_x ^2
            \end{array}
        \right] 
        C
        -
        C^{\top} 
        \left[ 
            \begin{array}{cc}
                0 &  \\
                  & \S_y ^2
            \end{array}
        \right] 
        C 
        u
        -\lbd u
        =
        F.
    \end{equation*}
    A further transformation by using $\tilde{u} := CuC^{\top}$ and $\tilde{F} := CFC^{\top}$ results with
    \begin{equation*}
        -\tilde{u} 
        \left[ 
            \begin{array}{cc}
                0 &  \\
                  & \S_x ^2
            \end{array}
        \right]
        -
        \left[ 
            \begin{array}{cc}
                0 &  \\
                  & \S_y ^2
            \end{array}
        \right] 
        \tilde{u}
        -\lbd \tilde{u}
        =
        \tilde{F}.
    \end{equation*}
    The solution to the above equation can be obtained by a direct entrywise division due to linearity of the equation and the non-singularity of the coefficient matrix (non zero $\lbd$), giving 
    \begin{equation*}
        \tilde{u} = \tilde{F} ./ M,
    \end{equation*}
    where $M$ is the $N \times N$ coefficient matrix defined as
    \begin{equation*}
    \resizebox{0.99\columnwidth}{!}{$\displaystyle
        M = - \left[
            \begin{array}{ccccc}
                0                & \s_{1,x}^2                    & \cdots & \s_{N-2,x}^2                    & \s_{N-1,x}^2 \\
                \s_{1,y}^2   & \s_{1,x}^2 + \s_{1,y}^2   & \cdots & \s_{N-2,x}^2 + \s_{1,y}^2   & \s_{N-1,x}^2 + \s_{1,y}^2 \\
                \vdots           & \vdots                            & \ddots & \vdots                              & \vdots \\
                \s_{N-2,y}^2 & \s_{1,x}^2 + \s_{N-2,y}^2 & \cdots & \s_{N-2,x}^2 + \s_{N-2,y}^2 & \s_{N-1,x}^2 + \s_{N-2,y}^2 \\
                \s_{N-1,y}^2 & \s_{1,x}^2 + \s_{N-1,y}^2 & \cdots & \s_{N-2,x}^2 + \s_{N-1,y}^2 & \s_{N-1,x}^2 + \s_{N-1,y}^2 \\
            \end{array}
        \right] - \lbd
                            $}.
    \end{equation*}
    The solution to \cref{eqs:hlm} is thus
    \begin{equation*}
        u = C^{\top}((CFC^{\top})./M)C.
    \end{equation*}
    Note that the discrete cosine transformations are easy to implement in \textsc{Matlab} using commands \verb|dct2| for the forward transform ($C[\cdot]C^{\top}$) and \verb|idct2| for the inverse transform ($C^{\top}[\cdot]C$).        

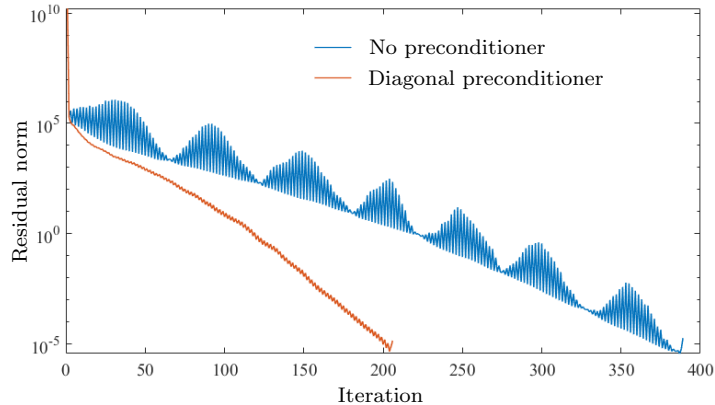
\begin{figure}[!htbp]
        \captionsetup[subfigure]{
                                justification=centering,
                                }
        \centering
        \begin{minipage}{1.0\textwidth}
            \centering
            \subfloat
                {
                \label{subfig:preconditioner}
                \begin{tikzpicture}
                \draw (0, 0) node[inner sep=0] {\includegraphics[width=0.7\linewidth,
                                angle=0]{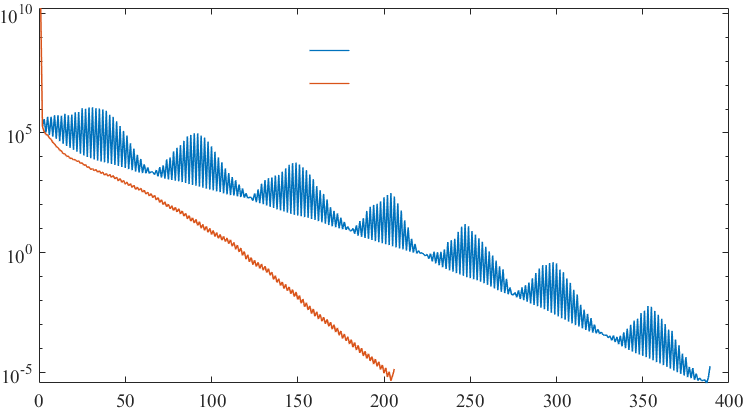}};
                \draw node[rotate=90] at (-4.7,0) {\footnotesize Residual norm};
                \draw (0.1, -2.7) node {\footnotesize Iteration};
                \draw (1.12, 1.9) node {\footnotesize No preconditioner};
                \draw (1.5, 1.5) node {\footnotesize Diagonal preconditioner};
                \end{tikzpicture}
                }
            \vspace{-0.02\linewidth}
        \end{minipage}
        \caption{Convergence of the preconditioned conjugate gradient method solving (\ref{eqs:sub-i-el})-(\ref{eqs:sub-i-el-bc}). Residual norm against iteration. The diagonal preconditioner is very effective in this case.}
        \label{fig:pcg}
    \end{figure}

\vspace{1.0em}
\noindent \textbf{The I-subproblem:} The fourth problem is to solve for $I$, freezing the other variables, in (\ref{eqs:opt_min}).
    \begin{equation}\label{eqs:sub_i}
        I^{*} 
        = 
        \mbox{arg }\min_I 
            \left \{ 
                \hth | I - I_{\S} |^2 
                - \Lbd_P \cdot \D I 
                + \cfrac{c_P}{2} | \mbf{P} - \D I |^2 
            \right \}.
    \end{equation}
    We will find a fast solver to this. The corresponding Euler-Lagrange equation is as the following.
    \begin{equation}\label{eqs:sub-i-el}
        \D \cdot (\D I^*)  - \cfrac{2 \hth}{c_P} I^*  =  \D \cdot \mbf{P} + \cfrac{1}{c_P} \D \cdot \Lbd_P - \cfrac{2 \hth}{c_P} I_{\S},
    \end{equation} with the Neumann boundary condition
    \begin{equation}\label{eqs:sub-i-el-bc}
        \D I^{*}  \cdot \mbf{\nu} =  (\mbf{P} + \frac{1}{c_P} \Lbd_P) \cdot \mbf{\nu},
    \end{equation} where the $\mbf{\nu}$ is the outward unit normal vector on the boundary of the domain. The above equation is not a standard modified Helmholtz equation because the coefficient ${2\hth}/{c_P}$ is a scalar function. We solve this equation with the conjugate gradient method, cf. e.g. \cite{NW2006}, in \textsc{Matlab} using the function \verb|pcg|. Using the diagonal preconditioner is the simplest, yet effective in this case, cf. \cref{fig:pcg}, a typical convergence result is shown.
        
\section{Numerical experiments}
\label{sec:numt}

In this section we present our experiments, on the choice of our vector matching term and the regularizers in the proposed model, on the effectiveness of the model on simple geometries, on its being able to accurately represent structures like edges and geometries, and finally on its effectiveness on real 3D maps with very few level lines comparing it with the state of the art model. 

We use the augmented Lagrangian method of \cref{sec:alg} for the solution. The parameters of the augmented Lagrangian method are set as $c_Q = 1$, $c_P = 1$ and $c_E = 1$ throughout our experiment, unless otherwise stated. The information needed are level curves or level lines, often provided as point clouds. In practice, these points are normally oriented, and is therefore easy to construct level lines from them. In case they are not oriented, we need to determine their orientations. To do that we simply connect the nearest given points with same level value to form the the level lines. In some cases the nearest points cannot be identified, for which, we first run the isotropic model ($\a=0$ in (\ref{eqs:opt_normal})), corresponding to the model without vector matching. Once the isotropic surface is obtained, the \textsc{Matlab} function called \verb|contour| is used to get the level lines which are then used as guidelines to find the order of the given points. A threshold is used to determine the connectivity between the points, once exceeded the level lines are then considered disconnected.

\subsubsection*{On the vector matching term}
\begin{figure}[!htbp]
        \captionsetup[subfigure]{
                                justification=centering,
                                labelformat=empty,
                                }
        \centering
        \begin{minipage}{1.0\textwidth}
            \centering
            \subfloat[][$\a = 0$\\]
                {
                \includegraphics[width=0.165\linewidth,
                                height=0.165\linewidth,
                                angle=0]{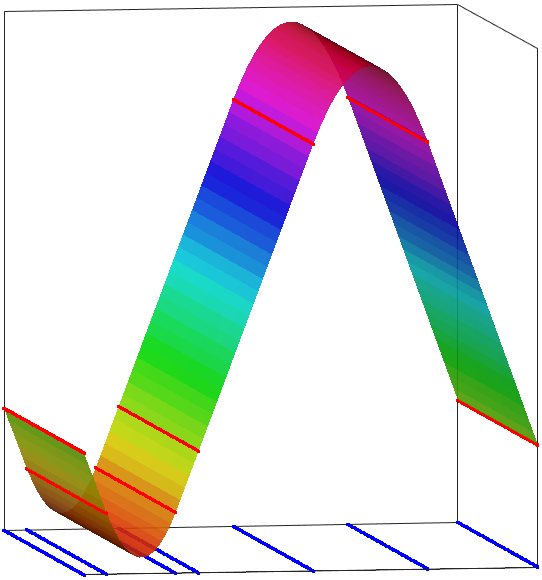}
                \label{subfig:pln_t0}
                }
            \hspace{0.01\linewidth}
            \subfloat[][$\a = 10^5$\\]
                {
                \includegraphics[width=0.165\linewidth,
                                height=0.165\linewidth,
                                angle=0]{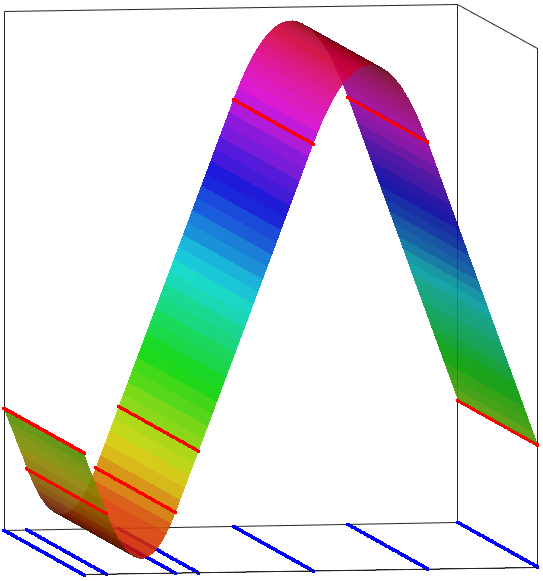}
                \label{subfig:pln_t100000}
                }
            \hspace{0.01\linewidth}
            \subfloat[][$\a = 10^2$ \\ $\angle (\mbf{v}_{\G}, \boldsymbol{\tau}_{\G}) = \ang{5}$]
                {
                \includegraphics[width=0.165\linewidth,
                                height=0.165\linewidth,
                                angle=0]{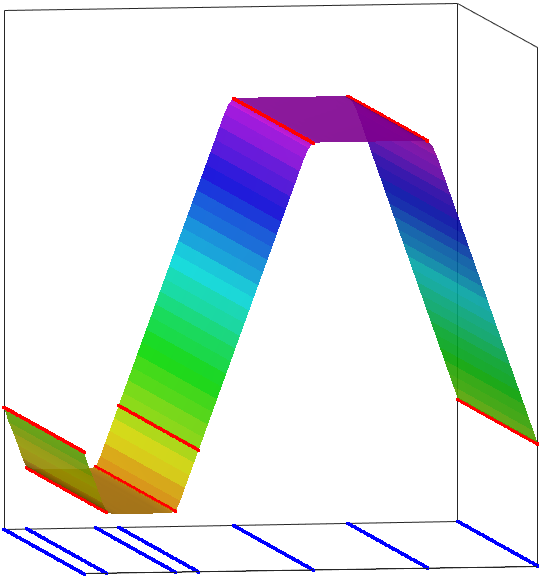}
                \label{subfig:pln_t100}
                }
            \hspace{0.01\linewidth}
            \subfloat[][$\a = 1$ \\ $\mbf{v}_{\G} \perp \boldsymbol{\tau}_{\G}$]
                {
                \includegraphics[width=0.165\linewidth,
                                height=0.165\linewidth,
                                angle=0]{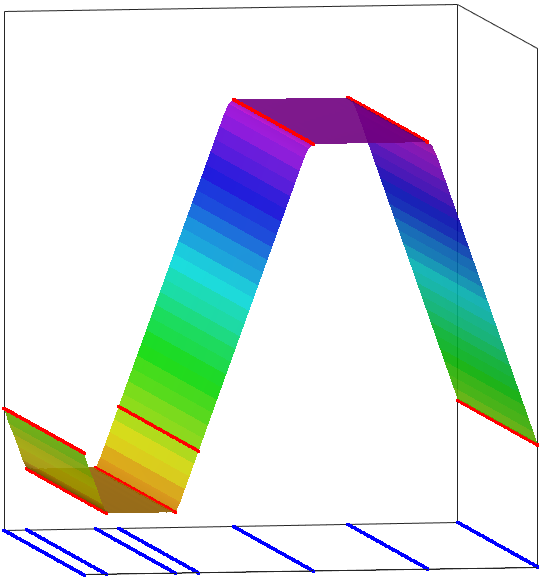}
                \label{subfig:pln_t1}
                }
            \vspace{-0.01\linewidth}
        \end{minipage}
        \caption{Illustrating the use of tangent vector matching term (\ref{eqs:tangent}) in the model (\ref{eqs:opt}) to get (\ref{eqs:opt_tangent}). Surfaces are reconstructed from the level lines (the blue lines on the floor, or the red lines on the surface) with the vector $\mbf{v}_{\G}$ prescribed only on the level lines with the largest and the smallest height values. Here, only the second order regularization is used. $\boldsymbol{\tau}_{\G}$ is the tangent vector of the given level lines. The first two figures correspond to the case $\D I^* \perp \mbf{v}_{\G}$ (no matter how big the value of $\a$ is, the surface is the same isotropic surface) and the last two figures correspond to the case $\D I^*=0$ (in which case the surface is a flattened structure).}
        \label{fig:v_tangent}
    \end{figure}

\begin{figure}[!htbp]
        \captionsetup[subfigure]{
                                justification=centering,
                                labelformat=empty,
                                }
        \centering
        \begin{minipage}{1.0\textwidth}
            \centering
            \subfloat[][$\a = 1.5$\\]
                {
                \includegraphics[width=0.165\linewidth,
                                height=0.165\linewidth,
                                angle=0]{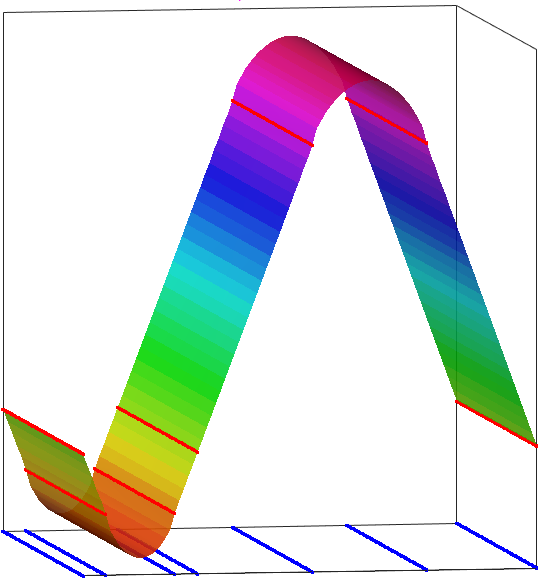}
                \label{subfig:pln_n1p5}
                }
            \hspace{0.01\linewidth}
            \subfloat[][$\a = 2$\\]
                {
                \includegraphics[width=0.165\linewidth,
                                height=0.165\linewidth,
                                angle=0]{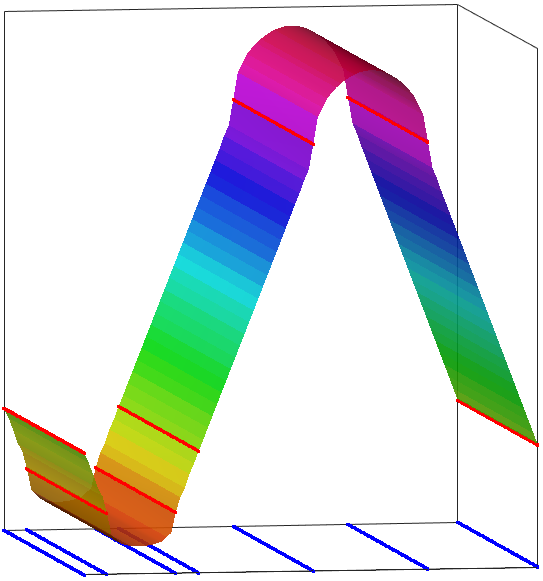}
                \label{subfig:pln_n2}
                }
            \hspace{0.01\linewidth}
            \subfloat[][ $\a = -0.1$\\]
                {
                \includegraphics[width=0.165\linewidth,
                                height=0.165\linewidth,
                                angle=0]{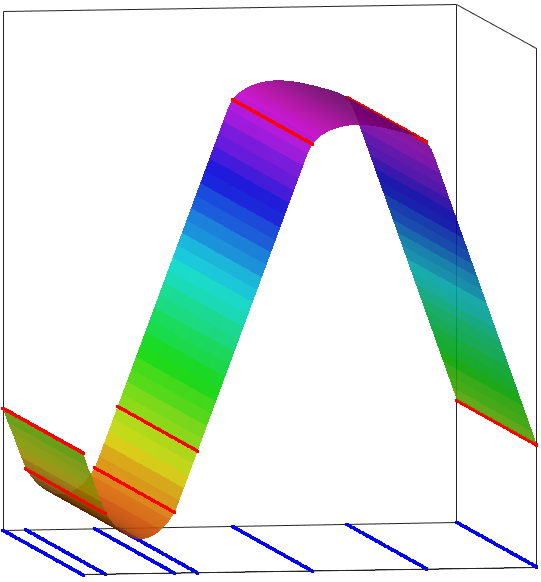}
                \label{subfig:pln_nn0p1}
                }
            \hspace{0.01\linewidth}
            \subfloat[][$\a = -0.2$\\]
                {
                \includegraphics[width=0.165\linewidth,
                                height=0.165\linewidth,
                                angle=0]{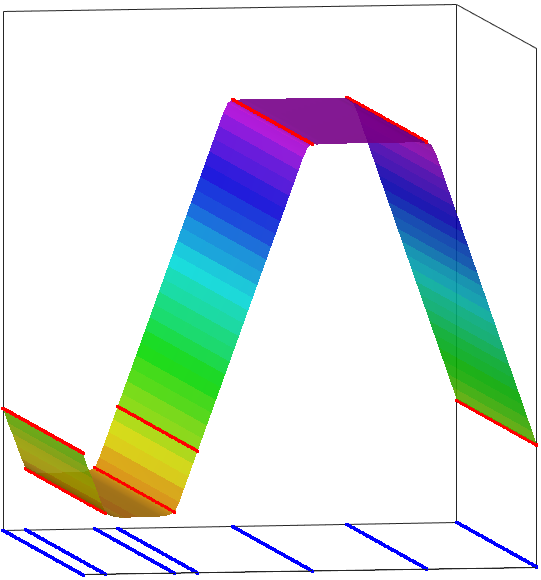}
                \label{subfig:pln_nn0p2}
                }
            \vspace{-0.01\linewidth}
        \end{minipage}
        \caption{Illustrating the use of normal vector matching (\ref{eqs:normal}) in the model (\ref{eqs:opt}), giving (\ref{eqs:opt_normal}). Surfaces are reconstructed from given level lines (the blue lines on the floor, or the red lines on the surface), with unit normal vectors prescribed only on the level lines with the largest and the smallest height values. Here, only the second order regularizer is used. The different shapes are reconstructed using different values of $\a$.}
        \label{fig:v_normal}
    \end{figure}

    With this experiment we justify our choice of vector matching term $V(\D I,\mbf{v}_{\G})$ in the minimization (\ref{eqs:opt}), from the two options (\ref{eqs:tangent}) and (\ref{eqs:normal}). To see the difference, we choose the following simple example: level lines (with height value) parallel to the $y$-axis, cf. \cref{fig:v_tangent} and \cref{fig:v_normal} (blue lines on the floor or red lines on the surface). Vector $\mbf{v}_{\G}$ is defined only on the lowest and highest level lines. The regularization and data fidelity parameters are kept the same as $g=1$, $h=0$ and $\te=10^5$ throughout this experiment.

    We consider the tangent vector matching first, model (\ref{eqs:opt_tangent}), this is illustrated in \cref{fig:v_tangent}. As mentioned earlier, there are two solutions, either $\D I^* = \mbf{0}$, which is reflected in the last two sub-figures, or $\D I^* \perp \mbf{v}_{\G}$ with $\D I^* \neq \mbf{0}$, which is reflected in the first two sub-figures. There are basically two shapes that can be obtained through the model (\ref{eqs:opt_tangent}), cf. \cref{fig:v_tangent}, flat structures being created between the level lines where vectors are given, and parabolic structures otherwise. $\mbf{v}_{\G}$ is chosen so that it is not parallel to tangent $\boldsymbol{\tau}_{\G}$ (to the level lines) in the last two sub-figures.

\begin{figure}[!htbp]
        \captionsetup[subfigure]{
                                justification=centering,
                                labelformat=empty,
                                }
        \centering
        \begin{minipage}{1.0\textwidth}
            \centering
            \subfloat[][$\a = - 1$\\]
                {
                \includegraphics[width=0.165\linewidth,
                                height=0.165\linewidth,
                                angle=0]{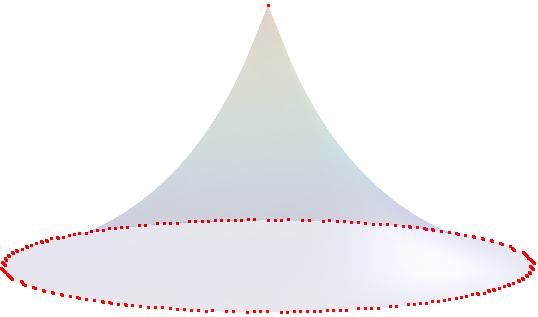}
                }
            \hspace{0.01\linewidth}
            \subfloat[][$\a = 0$\\]
                {
                \includegraphics[width=0.165\linewidth,
                                height=0.165\linewidth,
                                angle=0]{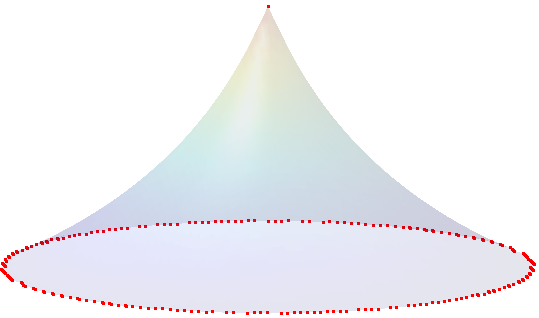}
                }
            \hspace{0.01\linewidth}
            \subfloat[][ $\a = 1$\\]
                {
                \includegraphics[width=0.165\linewidth,
                                height=0.165\linewidth,
                                angle=0]{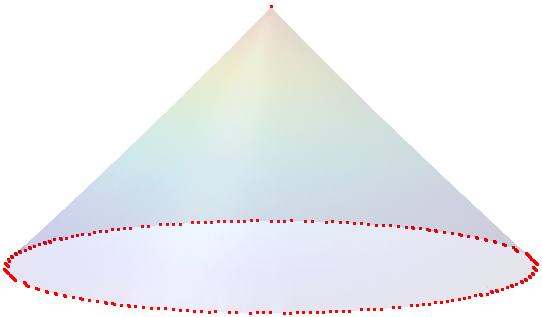}
                }
            \hspace{0.01\linewidth}
            \subfloat[][$\a = 2$\\]
                {
                \includegraphics[width=0.165\linewidth,
                                height=0.165\linewidth,
                                angle=0]{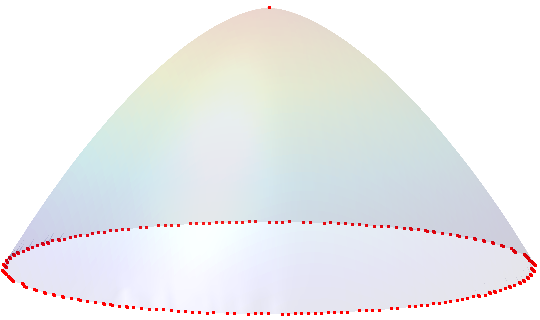}
                }
            \vspace{-0.01\linewidth}
        \end{minipage}
        \caption{Illustrating the use of normal vector matching (\ref{eqs:normal}) in our model (\ref{eqs:opt}), giving (\ref{eqs:opt_normal}). Surfaces are reconstructed from one level line (the red line on the floor) and one height in the center (the red point on the top), plus unit normal vectors prescribed only on the level line with the smallest height value. Here, only the second order regularizer is used. Cones with different shapes, from concave to convex, are reconstructed using different values of $\a$.}
        \label{fig:v_normal_2}
    \end{figure}

    We now consider the normal vector matching term, giving as the model (\ref{eqs:opt_normal}). The results are shown in \cref{fig:v_normal}. As shown in the figure, varying the parameter $\a$, we get different shapes, from parabolic to flattened. Note that, based on the normal vector matching, we are able to generate flat structures even when the first order regularizer has been inactive (corresponding to $h = 0$), see the last sub-figure. The value of the parameter $\a$ for the vector matching term is set equal to $1.5$, $2$, $-0.1$ and $-0.2$ respectively in \cref{fig:v_normal}. A second example showing similar results is presented in \cref{fig:v_normal_2}, where cones of different shapes, from concave to convex, have been obtained by varying the parameter $\alpha$ ($\alpha = -1, 0, 1,$ and $2$).

\subsubsection*{On the regularizer terms}
\begin{figure}[htb]
        \captionsetup[subfigure]{
                                justification=centering,
                                labelformat=empty,
                                }
        \centering
        \begin{minipage}{1.0\textwidth}
            \centering
            \subfloat
                {
                \includegraphics[width=0.165\linewidth,
                                height=0.165\linewidth,
                                angle=0]{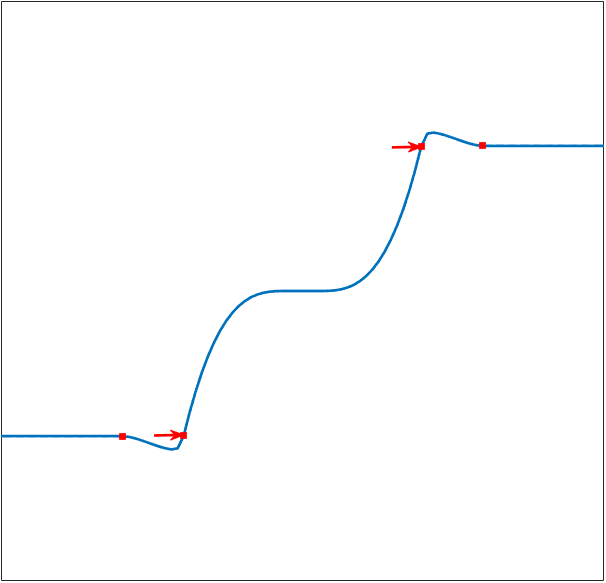}
                \label{subfig:ord1-1}
                }
            \hspace{0.005\linewidth}
            \subfloat
                {
                \includegraphics[width=0.165\linewidth,
                                height=0.165\linewidth,
                                angle=0]{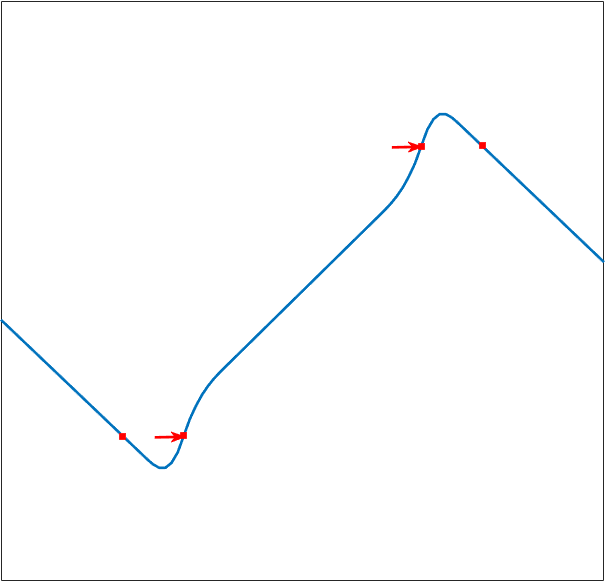}
                \label{subfig:ord2-1}
                }
            \hspace{0.005\linewidth}
            \subfloat
                {
                \includegraphics[width=0.165\linewidth,
                                height=0.165\linewidth,
                                angle=0]{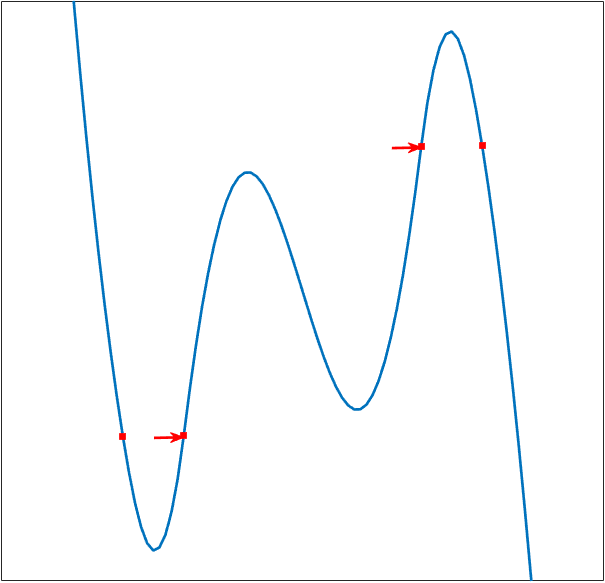}
                \label{subfig:ord3-1}
                }
            \hspace{0.005\linewidth}
            \subfloat
                {
                \includegraphics[width=0.165\linewidth,
                                height=0.165\linewidth,
                                angle=0]{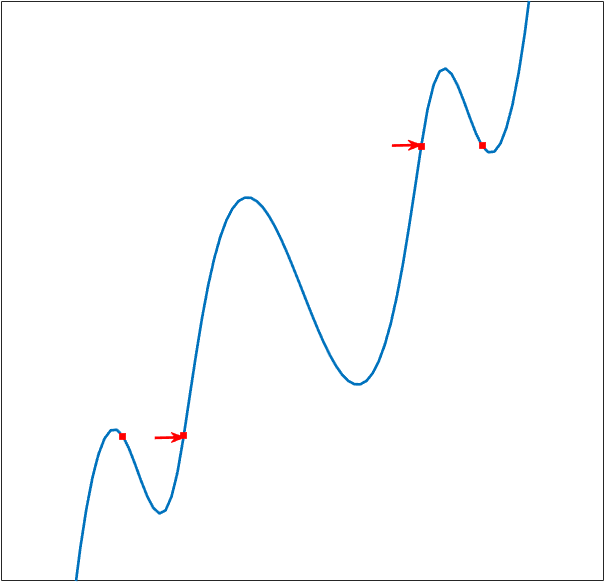}
                \label{subfig:ord4-1}
                }
            \hspace{0.005\linewidth}
            \subfloat
                {
                \includegraphics[width=0.165\linewidth,
                                height=0.165\linewidth,
                                angle=0]{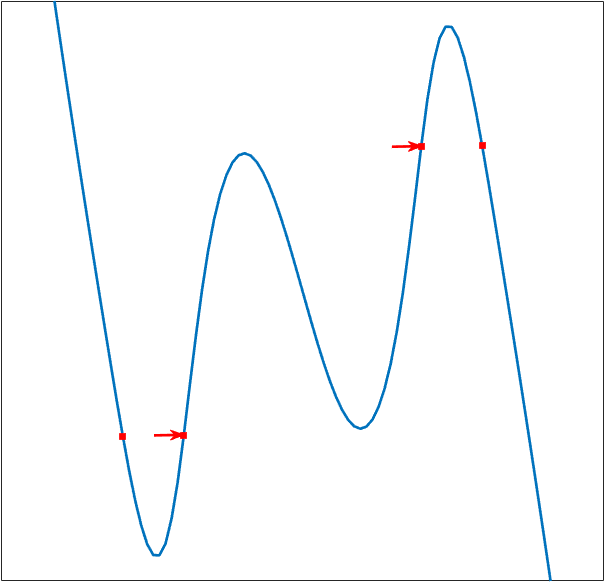}
                \label{subfig:ord5-1}
                }
            \vspace{-0.005\linewidth}
        \end{minipage}
        \begin{minipage}{1.0\textwidth}
            \centering
            \subfloat
                {
                \includegraphics[width=0.165\linewidth,
                                height=0.165\linewidth,
                                angle=0]{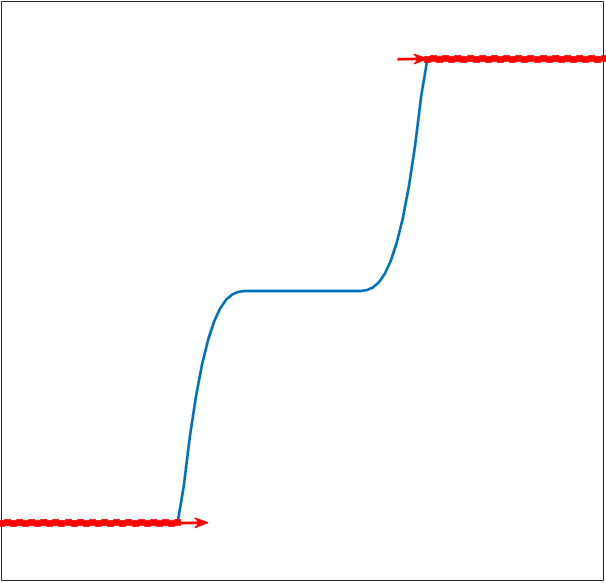}
                \label{subfig:ord1-2}
                }
            \hspace{0.005\linewidth}
            \subfloat
                {
                \includegraphics[width=0.165\linewidth,
                                height=0.165\linewidth,
                                angle=0]{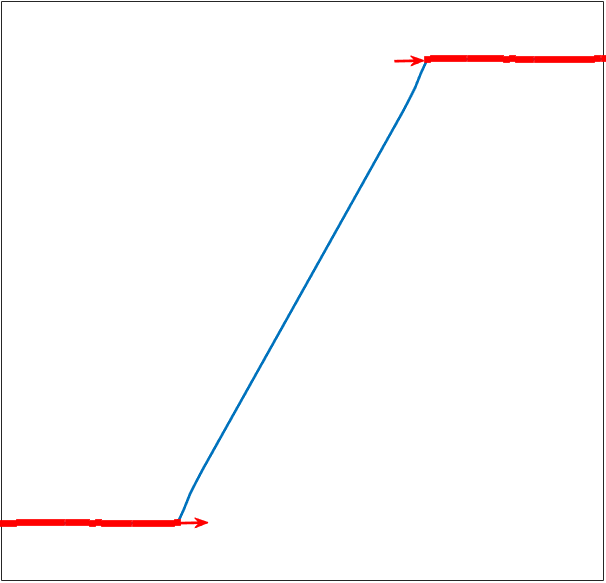}
                \label{subfig:ord2-2}
                }
            \hspace{0.005\linewidth}
            \subfloat
                {
                \includegraphics[width=0.165\linewidth,
                                height=0.165\linewidth,
                                angle=0]{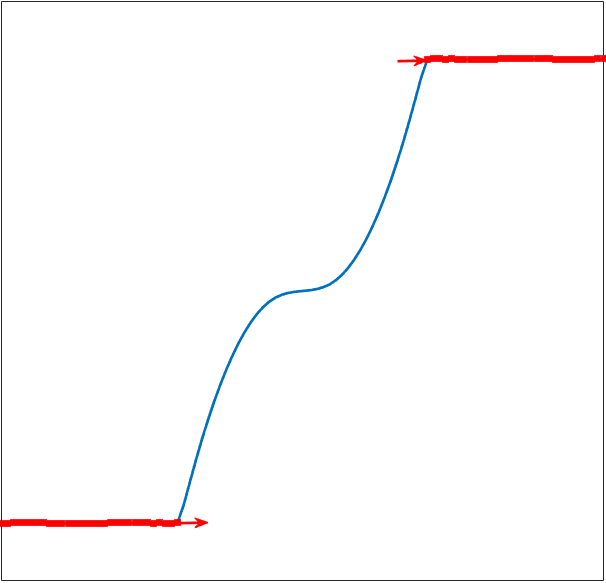}
                \label{subfig:ord3-2}
                }
            \hspace{0.005\linewidth}
            \subfloat
                {
                \includegraphics[width=0.165\linewidth,
                                height=0.165\linewidth,
                                angle=0]{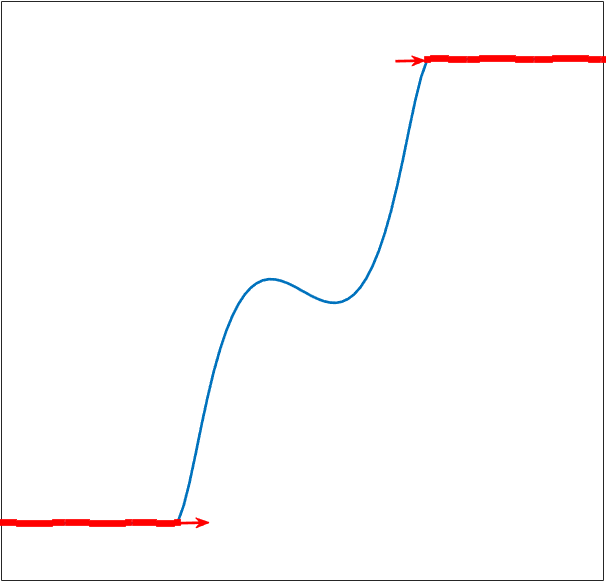}
                \label{subfig:ord4-2}
                }
            \hspace{0.005\linewidth}
            \subfloat
                {
                \includegraphics[width=0.165\linewidth,
                                height=0.165\linewidth,
                                angle=0]{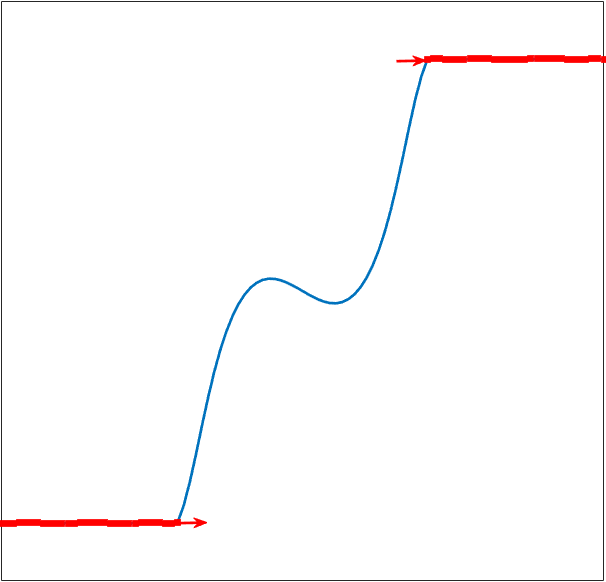}
                \label{subfig:ord5-2}
                }
            \vspace{-0.005\linewidth}
        \end{minipage}
        \clearsubcaptcounter
        \begin{minipage}{1.0\textwidth}
            \centering
            \subfloat[][$| \D^1 I |_F$\\]
                {
                \includegraphics[width=0.165\linewidth,
                                height=0.165\linewidth,
                                angle=0]{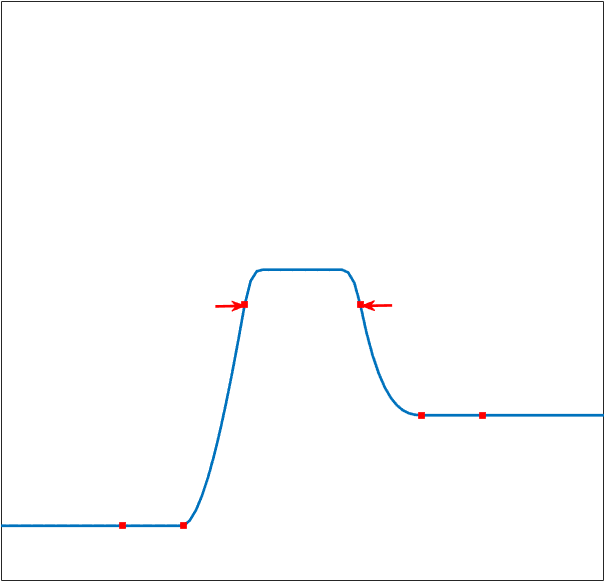}
                \label{subfig:ord1-4}
                }
            \hspace{0.005\linewidth}
            \subfloat[][$| \D^2 I |_F$\\]
                {
                \includegraphics[width=0.165\linewidth,
                                height=0.165\linewidth,
                                angle=0]{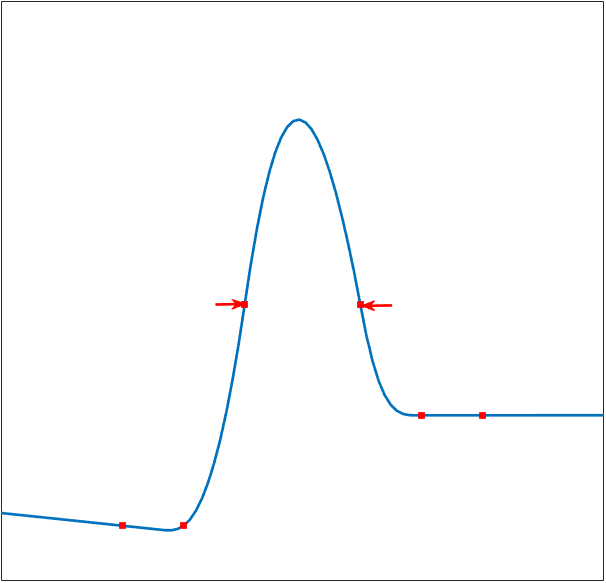}
                \label{subfig:ord2-4}
                }
            \hspace{0.005\linewidth}
            \subfloat[][$| \D^3 I |_F$\\]
                {
                \includegraphics[width=0.165\linewidth,
                                height=0.165\linewidth,
                                angle=0]{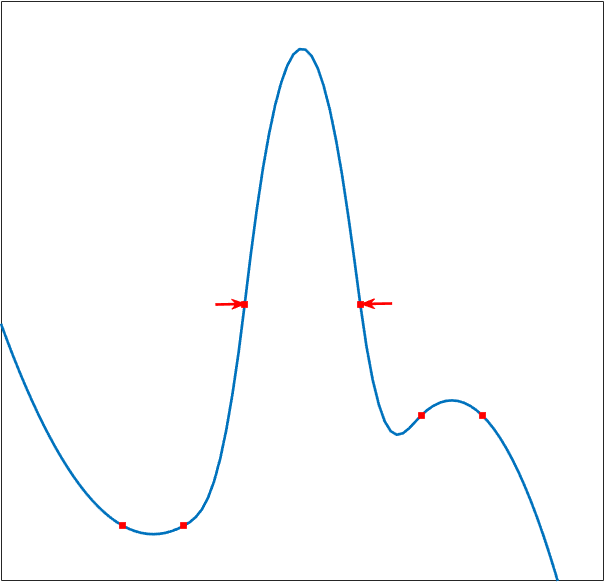}
                \label{subfig:ord3-4}
                }
            \hspace{0.005\linewidth}
            \subfloat[][$| \D^4 I |_F$\\]
                {
                \includegraphics[width=0.165\linewidth,
                                height=0.165\linewidth,
                                angle=0]{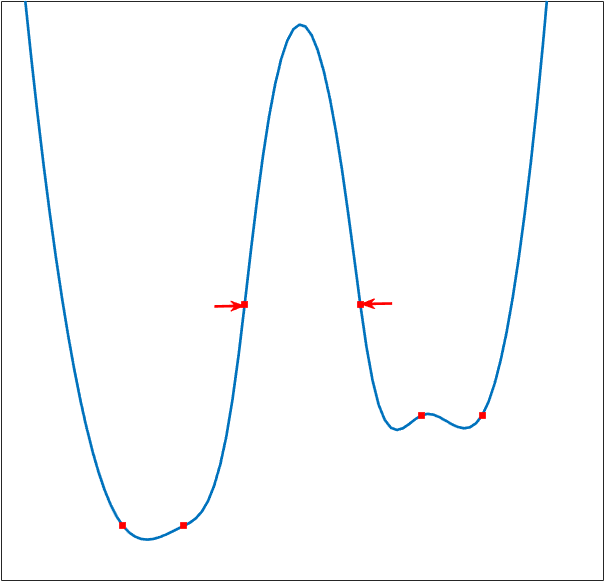}
                \label{subfig:ord4-4}
                }
            \hspace{0.005\linewidth}
            \subfloat[][$| \D^5 I |_F$\\]
                {
                \includegraphics[width=0.165\linewidth,
                                height=0.165\linewidth,
                                angle=0]{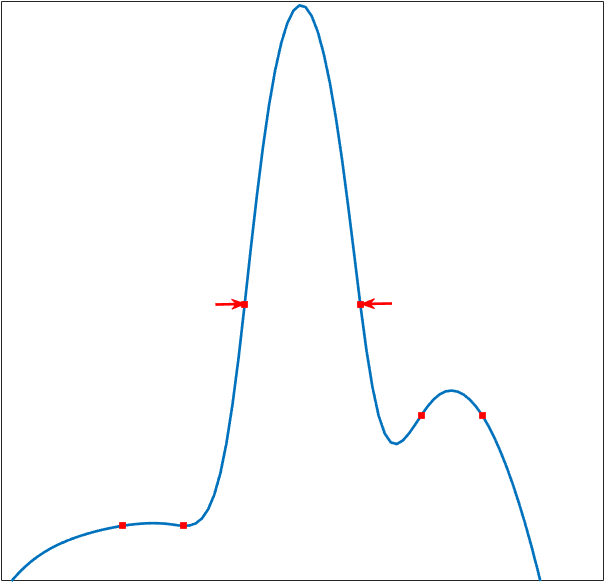}
                \label{subfig:ord5-4}
                }
            \vspace{-0.005\linewidth}
        \end{minipage}
        \caption{Showing The effect of using different order of regularizer in the reconstruction in 1D. The profiles are reconstructed using (\ref{eqs:opt_k_ord}); sparse data given are the height values (in red), and gradient direction using arrows (in blue). Each row corresponds to a specific test, and the columns corresponds to the orders from 1 to 5 of the regularizer. The first order results in flattened structure while the second order favours slopes in the structure. The higher order regularizer introduces oscillation in the solution, and makes it unbounded at the boundary.}
        \label{fig:ord_tv_vct}
    \end{figure}

    The second experiment is to study the effect of different orders of the regularizer. To do this we use only one regularizer term in our model and fix the vector matching term. The following minimization, same as (\ref{eqs:opt_k_ord}), is considered
    \begin{equation*}
        \min_I g | \D^k I |_{F \ (\O)}  - \a \D I \cdot \mbf{v}_{\G \ (\G)} + \te | I - I_{\S} |_{2 \ (\S)}^2,
        \label{eqs:opt_1d_vct}
    \end{equation*}
    where the first term is the $k^{th}$ order total variation (TV), the second term is the normal vector matching term, and the third term is the data fidelity term. $g$, $\a$ and $\te$ are the scalar parameters. $\G$ and $\S$ are the set points where the unit normal vector and elevation/height value are given, respectively. For clarity, we apply the model to one space dimension (1D). The results are shown in \cref{fig:ord_tv_vct}, each column represents a fixed order from $k=1$ to $5$, and each row corresponds to a particular test. The height values are given at red points and the vectors are marked with red arrows. Note that, in 1D, the unit vector is the value one with a sign for the direction, i.e. $\pm 1$. In \cref{fig:ord_tv_vct}, the positive sign corresponds to the arrow pointing to the right and vice versa. The parameters are kept the same throughout this experiment as $g=1$ for the regularizer, $\a=1$ for the vector matching term, and $\te=10$ for the data fidelity.

    The experiment shows that orders higher than two, when combined with the normal vector matching, may induce oscillation in the shape of the reconstructed surface, as well as unbounded surface at the boundary. A combination of the first order and the second order regularizer is enough to reconstruct most shapes. This will see in the following experiments.

\subsubsection*{The proposed model on regular structures}
    The third experiment is to validate the proposed model (\ref{eqs:opt_normal}) where both regularizers and the normal vector matching are used. We study its effectiveness in accurately representing simple geometries like edges and corners.

\begin{figure}[!htbp]
        \captionsetup[subfigure]{
                                justification=centering,
                                labelformat=empty,
                                }
        \centering
        \begin{minipage}{1.0\textwidth}
            \centering
            \subfloat[][Ground truth\\]
                {
                \includegraphics[width=0.165\linewidth,
                                height=0.165\linewidth,
                                angle=0]{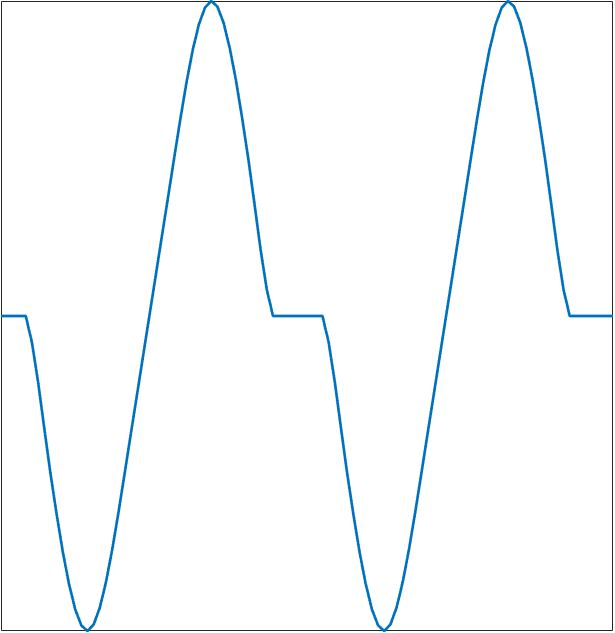}
                \label{subfig:sgn_gt}
                }
            \hspace{0.005\linewidth}
            \subfloat[][$| \D^2 I |_F$\\]
                {
                \includegraphics[width=0.165\linewidth,
                                height=0.165\linewidth,
                                angle=0]{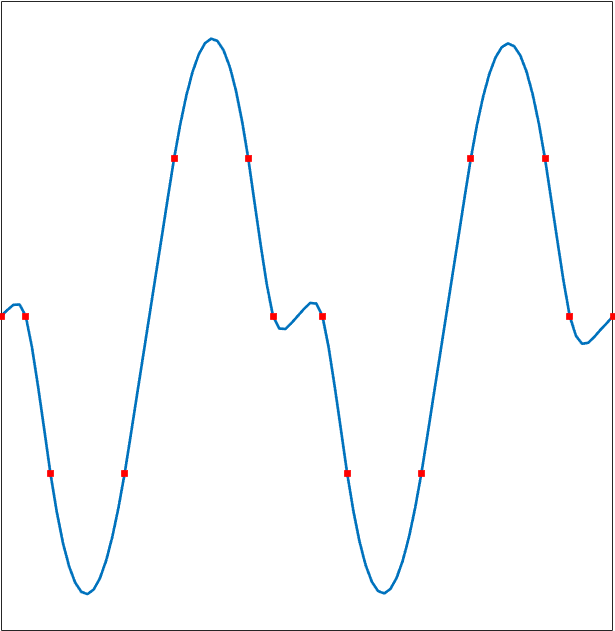}
                \label{subfig:sgn_allg}
                }
            \hspace{0.005\linewidth}
            \subfloat[][$| \D^1 I |_F$\\]
                {
                \includegraphics[width=0.165\linewidth,
                                height=0.165\linewidth,
                                angle=0]{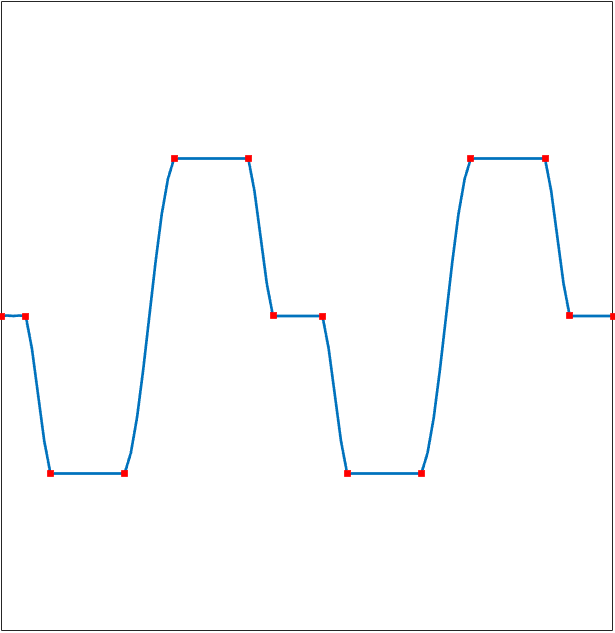}
                \label{subfig:sgn_allh}
                }
            \hspace{0.005\linewidth}
            \subfloat[][$| \D^1 I |_F + | \D^2 I |_F$\\]
                {
                \includegraphics[width=0.165\linewidth,
                                height=0.165\linewidth,
                                angle=0]{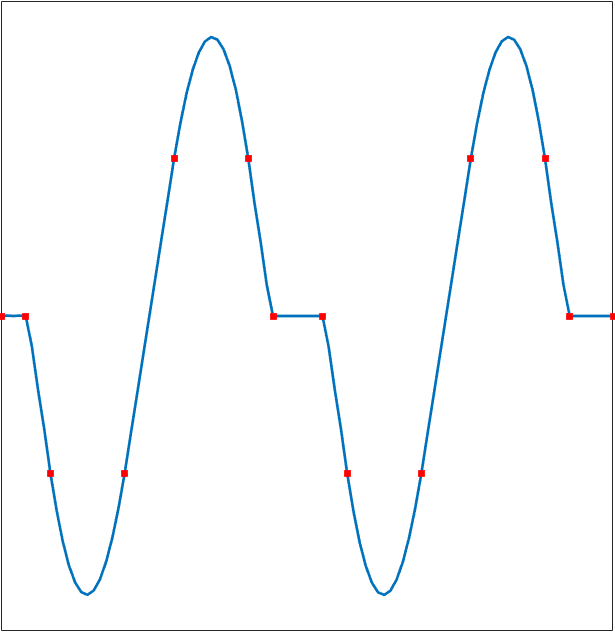}
                \label{subfig:sgn_gh}
                }
            \hspace{0.005\linewidth}
            \subfloat[][Proposed model\\]
                {
                \includegraphics[width=0.165\linewidth,
                                height=0.165\linewidth,
                                angle=0]{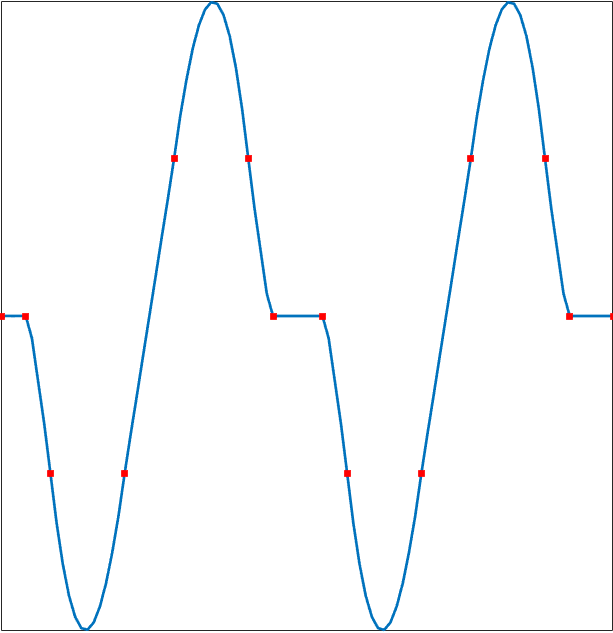}
                \label{subfig:sgn_vct}
                }
            \vspace{-0.01\linewidth}
        \end{minipage}
        \caption{Showing the effect of using vector matching in combination with both first order and second order total variation regularizer in the proposed model (\ref{eqs:opt_normal}), as compared to the second, third and fourth figures from left where the vector matching has not been used.}
        \label{fig:tvs_vct}
    \end{figure}

    We first test the proposed model on a 1D signal with the mixed shape of sine and rectangular waves, cf. \cref{fig:tvs_vct}. Only a few isolated points (red points) with height values and unit normal vectors ($\pm 1$) are given. The normal vector matching is only used in the last sub-figure, with the parameter $\a = 1.5$. Throughout this test, the parameter $h = 1.5$ for the first order total variation (TV), i.e., $| \D^1 I |_F$, and the parameter $g = 4$ for the second order TV, i.e., $| \D^2 I |_F$. 

    The result as shown in \cref{fig:tvs_vct}, demonstrates that the proposed model (\ref{eqs:opt_normal}) is able to accurately recover the shape with corners. While the first order TV favors the staircase structures and the second TV preserves sloped structures, the combination of the first order and the second order TVs without vector matching reconstructs both flattened and sloped structures quite well. However it still fails to be accurate on the top and in valleys.

\begin{figure}[!htbp]
        \captionsetup[subfigure]{
                                justification=centering,
                                labelformat=empty,
                                }
        \centering
        \begin{minipage}{1.0\textwidth}
            \centering
            \subfloat[][Ground truth\\]
                {
                \includegraphics[width=0.165\linewidth,
                                height=0.165\linewidth,
                                angle=0]{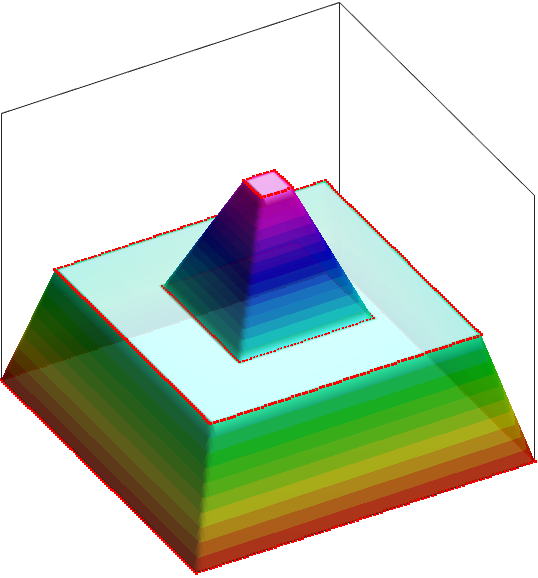}
                \label{subfig:stps_gt}
                }
            \hspace{0.01\linewidth}
            \subfloat[][$| \D^2 I |_F$\\]
                {
                \includegraphics[width=0.165\linewidth,
                                height=0.165\linewidth,
                                angle=0]{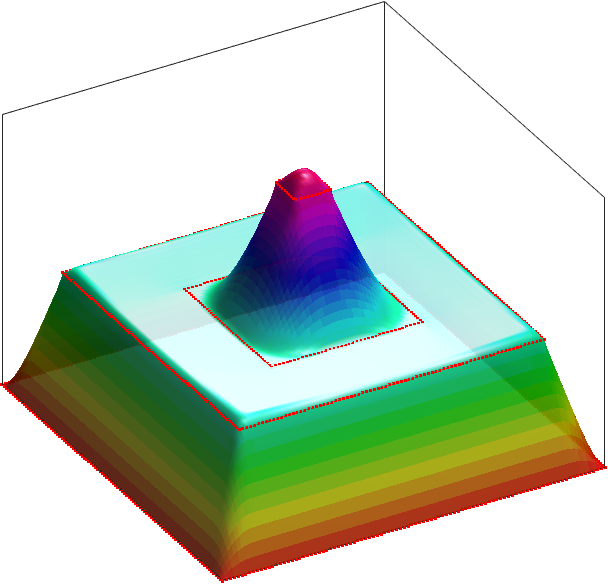}
                \label{subfig:stps_allg}
                }
            \hspace{0.01\linewidth}
            \subfloat[][$| \D^1 I |_F$\\]
                {
               \includegraphics[width=0.165\linewidth,
                                height=0.165\linewidth,
                                angle=0]{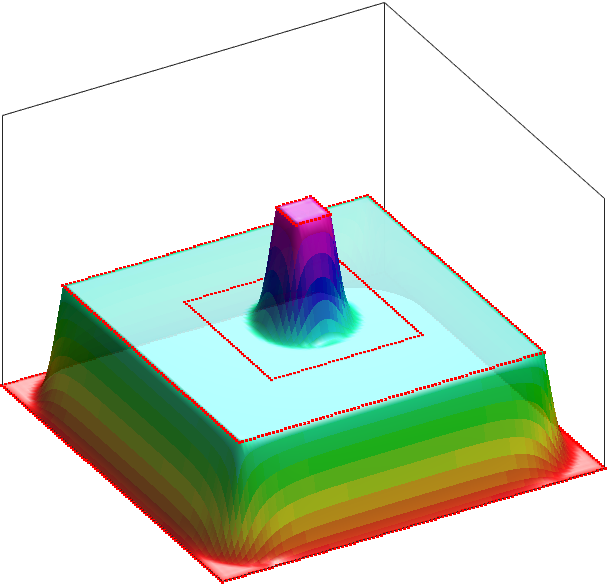}
                \label{subfig:stps_allh}
                }
            \hspace{0.01\linewidth}
            \subfloat[][$| \D^1 I |_F$ + $| \D^2 I |_F$\\]
                {
                \includegraphics[width=0.165\linewidth,
                                height=0.165\linewidth,
                                angle=0]{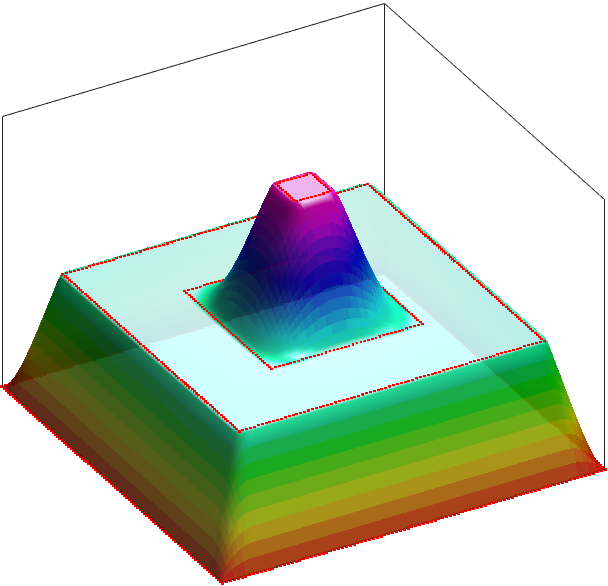}
                \label{subfig:stps_gh}
                }
            \hspace{0.01\linewidth}
            \subfloat[][Proposed model\\]
                {
                \includegraphics[width=0.165\linewidth,
                                height=0.165\linewidth,
                                angle=0]{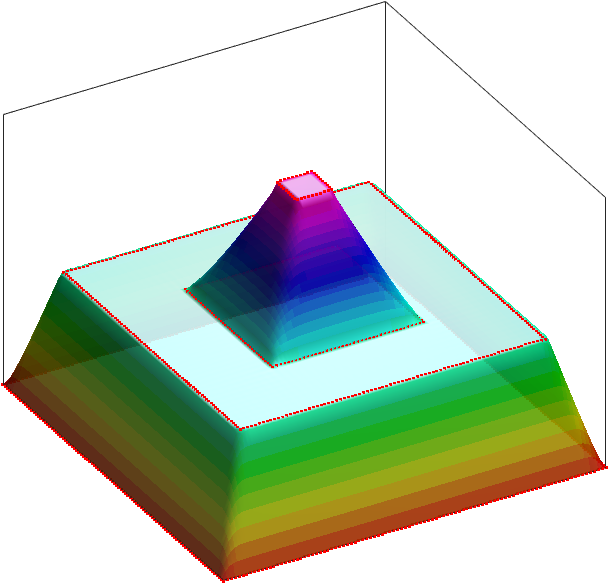}
                \label{subfig:stps_vct}
                }
            \vspace{-0.01\linewidth}
        \end{minipage}
        \caption{2D illustration showing the effect of using vector matching in combination with both first order and second order total variation regularizer in the proposed model (\ref{eqs:opt_normal}), in the last figure, as compared to the first, second and third figures from left where the vector matching has not been used.}
        \label{fig:stp_squares}
    \end{figure}

    Next, we apply our model in 2D, on a geometry that looks like the Maya pyramid, cf. \cref{fig:stp_squares}. There are three level lines available (red lines on the surface) with both height values and unit normal vectors, either given or extracted from the geometry of the underlying level lines. Note that, to accurately reconstruct the structures, we need to adjust the parameters of each term in our model (\ref{eqs:opt_normal}). The data fidelity parameter $\te$, in the experiment shown in \cref{fig:stp_squares}, is chosen to $\te=10^5$. The regularizer parameter $g = 10$ in the second sub-figure, while $h = 1$ in the third sub-figure. For sloped structures such as the walls of the base and the walls of the top, we simply use only the second order regularizer for the last two sub-figures, where $h=0$. In these two sub-figures, $g = 10$ and $g = 1$ for the walls of the base and the walls of the top, respectively. For flattened structures like the horizontal area of the pyramid, we use only the first order regularizer in the last two sub-figures, where $g=0$ and $h = 10$. The normal vector matching is only used in the last sub-figure, where $\a = 1$ along the level lines except the ones are intersection between the base and the top of the pyramid, where $\a = 11$. A higher value of the parameter $\a$ setting on this intersection is to get sharper edges for the top.

    Similar to the 1D test showing in \cref{fig:tvs_vct}, the first order TV favors the flat structures and the second TV preserves sloped structures. The combination of the first order and the second order TVs without vector matching reconstructs both flat and sloped structures but fails to accurately represent edges and corners. The proposed model (\ref{eqs:opt_normal}) reconstructs the Maya pyramid more accurately representing the edges and corners.

\begin{figure}[!htbp]
        \captionsetup[subfigure]{
                                justification=centering,
                                labelformat=empty,
                                }
        \centering
        \begin{minipage}{1.0\textwidth}
            \centering
            \subfloat[][$\a=2.5$\\]
                {
                \includegraphics[width=0.175\linewidth,
                                height=0.175\linewidth,
                                angle=0]{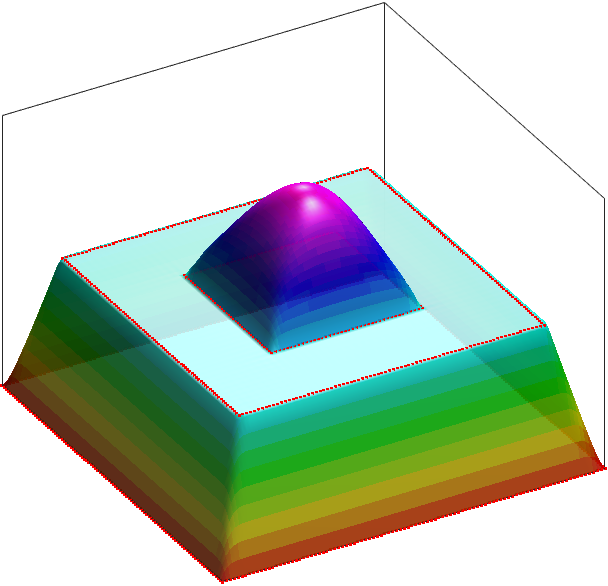}
                \label{subfig:vct_p2p5}
                }
            \hspace{0.01\linewidth}
            \subfloat[][$\a=2.0$\\]
                {
                \includegraphics[width=0.175\linewidth,
                                height=0.175\linewidth,
                                angle=0]{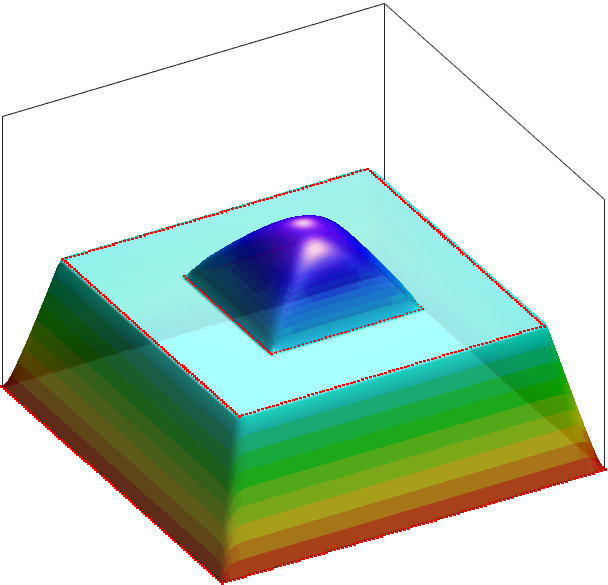}
                \label{subfig:vct_p2}
                }
            \hspace{0.01\linewidth}
            \subfloat[][$\a=0$\\]
                {
                \includegraphics[width=0.175\linewidth,
                                height=0.175\linewidth,
                                angle=0]{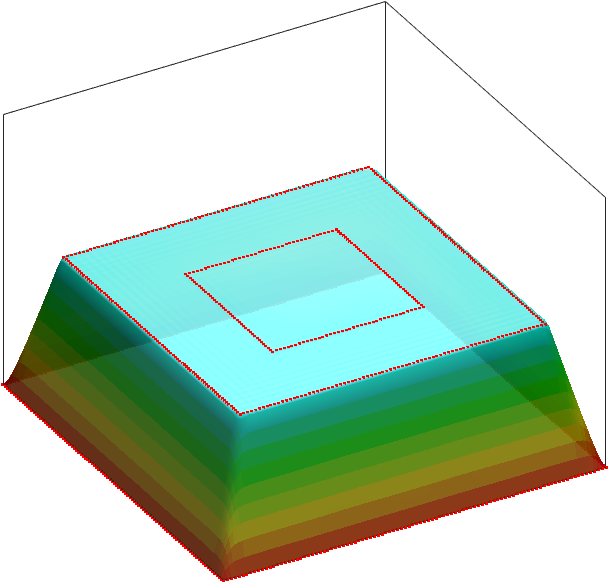}
                \label{subfig:vct_0}
                }
            \hspace{0.01\linewidth}
            \subfloat[][$\a=-2.0$\\]
                {
                \includegraphics[width=0.175\linewidth,
                                height=0.175\linewidth,
                                angle=0]{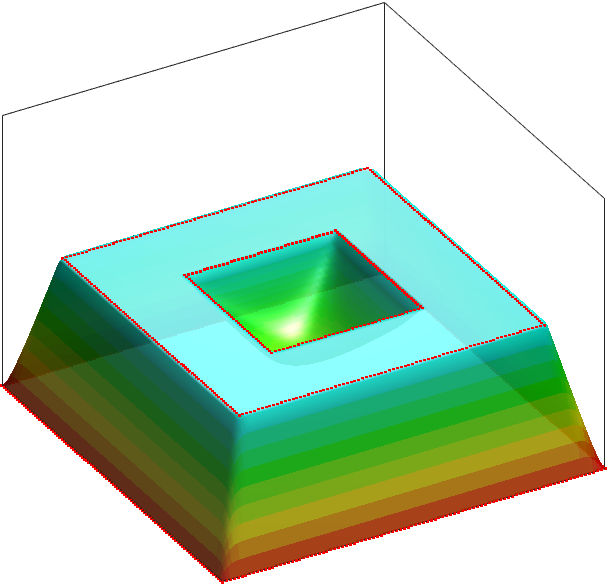}
                \label{subfig:vct_n2}
                }
            \vspace{-0.01\linewidth}
        \end{minipage}
        \caption{2D illustration showing how easily we can control the reconstructed surface through adjusting the vector matching parameter $\a$ in the proposed model (\ref{eqs:opt_normal}).}
        \label{fig:vct_alpha}
    \end{figure}

    The third test is to justify the capability of the proposed model (\ref{eqs:opt_normal}) in controlling of the reconstructed shape by adjusting the vector matching parameter $\a$. The results are shown in \cref{fig:vct_alpha}. Different to the previous test, there is no information given now at the center of the 2D domain, cf.\cref{fig:vct_alpha} and \cref{fig:stp_squares}. The parameters are kept the same as those in the last sub-figure of \cref{fig:stp_squares} except $\a$ on the innermost given level line (in red), cf. \cref{fig:stp_squares} where $\a = 2.5$, $\a = 2$, $\a = 0$ and $\a = -2$ for each test respectively. 
    The test shows that varying parameter $\a$ changes the shape, from a convex shape to a concave structure.

\begin{figure}[!htbp]
        \captionsetup[subfigure]{
                                justification=centering,
                                labelformat=empty,
                                }
        \centering
        \begin{minipage}{1.0\textwidth}
            \centering
            \subfloat[][1 contour\\]
                {
                \includegraphics[width=0.165\linewidth,
                                height=0.12\linewidth,
                                angle=0]{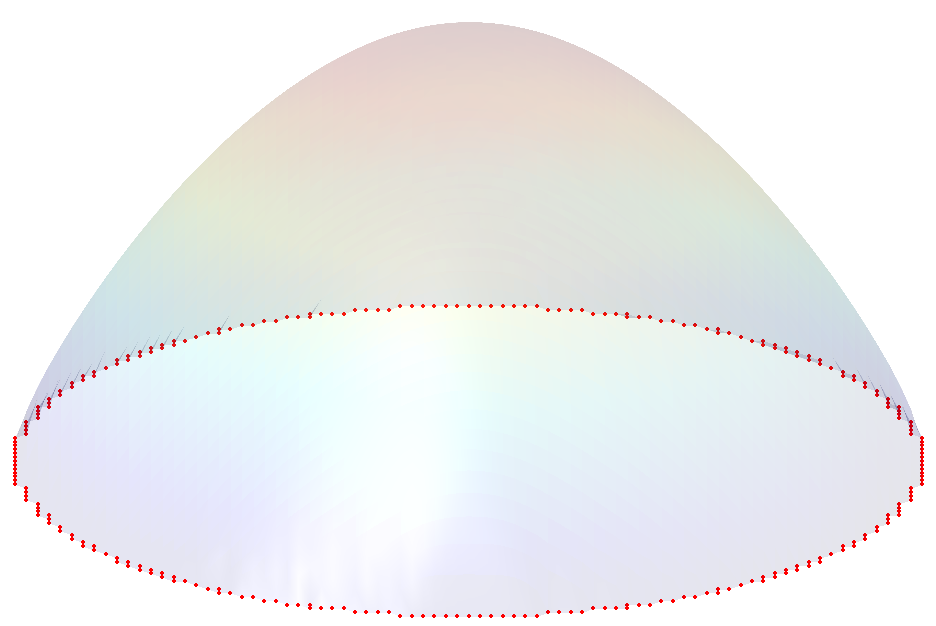}
                \label{subfig:ssph_1}
                }
            \hspace{0.005\linewidth}
            \subfloat[][2 contours\\]
                {
                \includegraphics[width=0.165\linewidth,
                                height=0.12\linewidth,
                                angle=0]{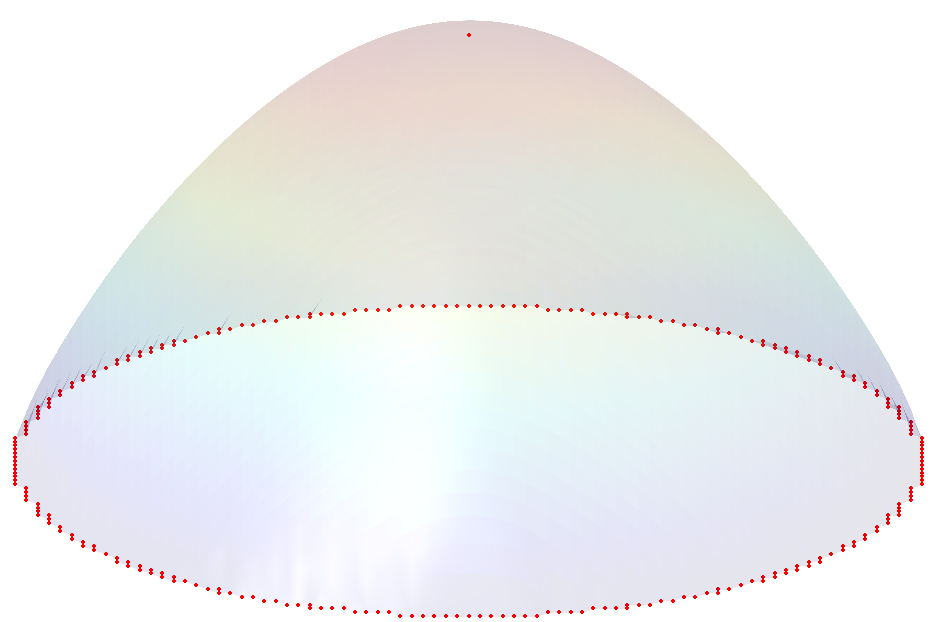}
                \label{subfig:ssph_2}
                }
            \hspace{0.005\linewidth}
            \subfloat[][4 contours\\]
                {
                \includegraphics[width=0.165\linewidth,
                                height=0.12\linewidth,
                                angle=0]{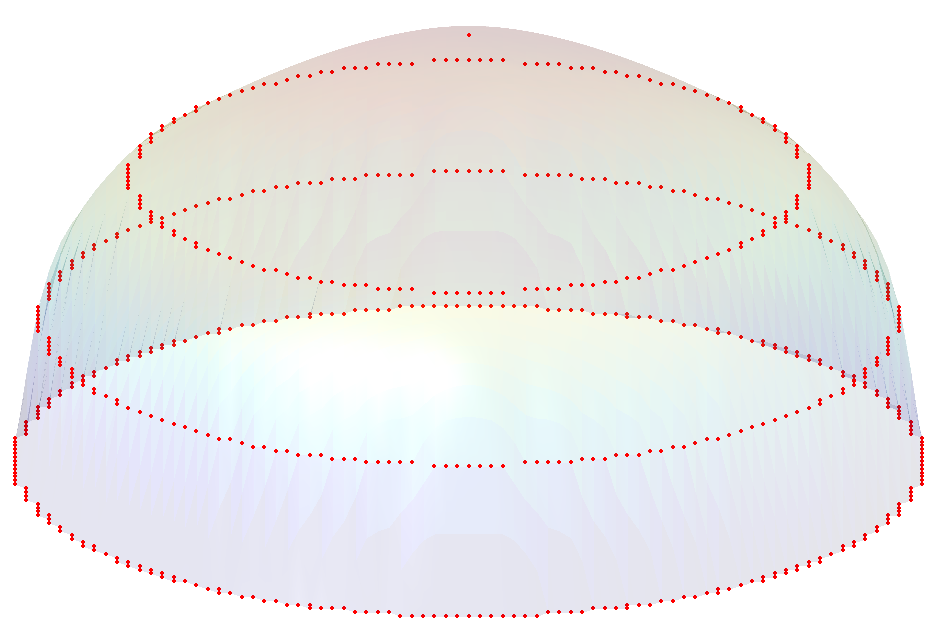}
                \label{subfig:ssph_4}
                }
            \hspace{0.005\linewidth}
            \subfloat[][8 contours\\]
                {
                \includegraphics[width=0.165\linewidth,
                                height=0.12\linewidth,
                                angle=0]{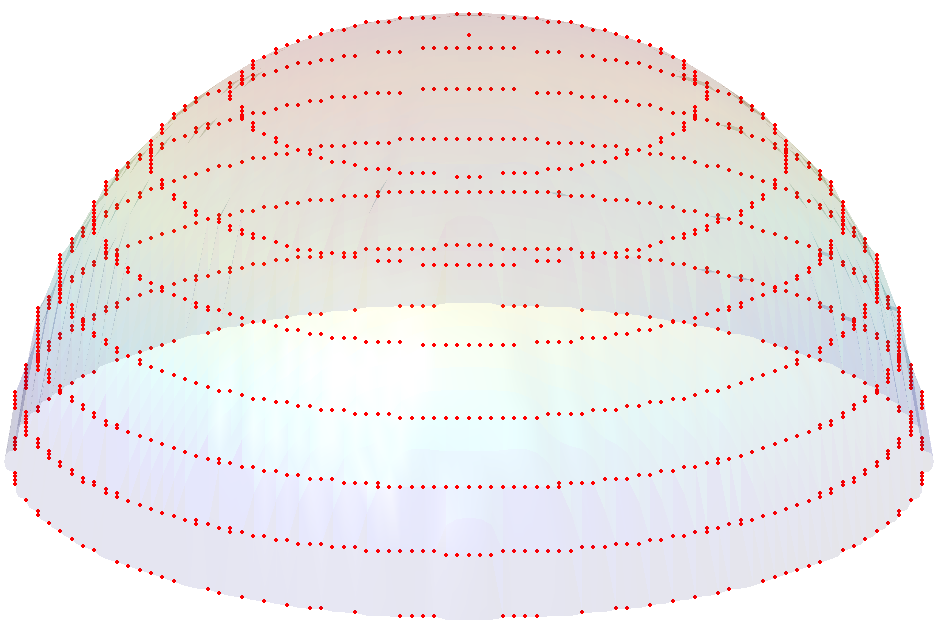}
                \label{subfig:ssph_8}
                }
            \vspace{-0.01\linewidth}
        \end{minipage}
        \caption{Showing the effect of level lines in the proposed model (\ref{eqs:opt_normal}), increasing the number of contours from left to right. With only one level line (the leftmost figure), the parameter $\a$ alone is enough to control the shape of the surface. The surface improves as the number of level lines increases.}
        \label{fig:ssph}
    \end{figure}

    The fourth test is on a semi-sphere. For each test we have different number of contours (the red lines on the surface), cf. \cref{fig:ssph}. Both the height value and the unit normal vectors are given on these contours. In this experiment, the regularizer parameters have been $g=1$ and $h=0$, and the data fidelity parameter has been $\te=10^5$. The vector matching parameter has been $\a=0.5$ along the contours, except the one with the largest radius, where $\a=2.85$.

    As shown here, even with only one single contour, as in the first sub-figure of \cref{fig:ssph}, the proposed model (\ref{eqs:opt_normal}) is still able to reconstruct the surface reasonably well. As we increase the number of contours, the reconstructed surface becomes more and more close to the perfect semi-sphere.

\subsubsection*{Real map I}
\begin{figure}[!htp]
        \captionsetup[subfigure]{
                                justification=centering,
                                labelformat=empty,
                                }
        \centering
        \begin{minipage}{1.0\textwidth}
            \centering
            \subfloat[][Three level lines]
                {
                \includegraphics[width=0.26\linewidth,
                                height=0.24\linewidth,
                                angle=0]{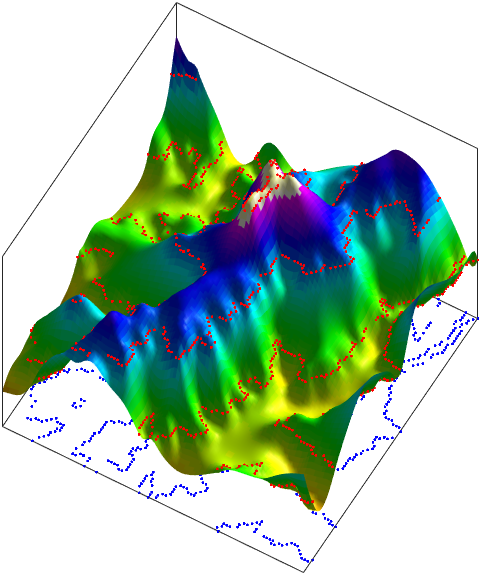}
                \label{subfig:mt_vct_3c}
                }
            \hspace{0.01\linewidth}
            \subfloat[][Seven level lines]
                {
                \includegraphics[width=0.26\linewidth,
                                height=0.24\linewidth,
                                angle=0]{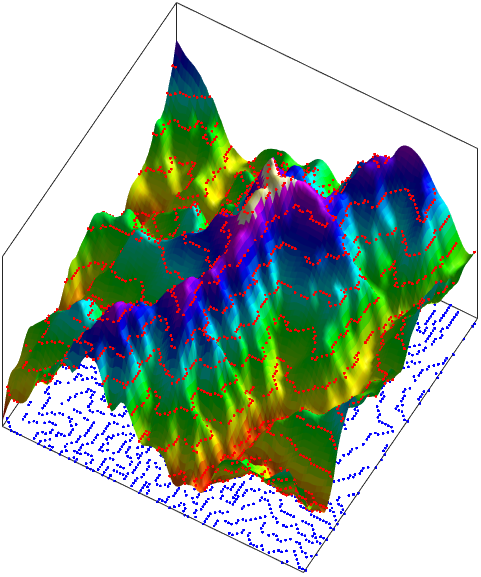}
                \label{subfig:mt_vct_7c}
                }
            \hspace{0.01\linewidth}
            \subfloat[][Ground truth]
                {
                \includegraphics[width=0.26\linewidth,
                                height=0.24\linewidth,
                                angle=0]{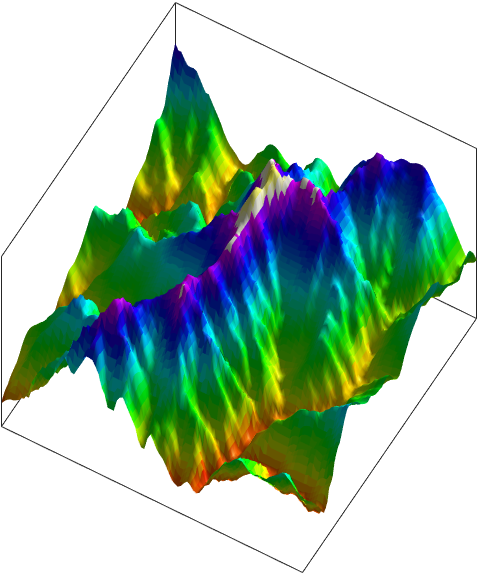}
                \label{subfig:mt_gt_3}
                }
            \vspace{-0.01\linewidth}
        \end{minipage}
        \caption{The real 3D map seen from the side, with increasing number of level lines from left to right. The surface reconstruction improves as the number of level lines increases.}
        \label{fig:mountain_sequence}
    \end{figure}

    In this experiment we apply our model to real maps with clear crease (kinks) and clear valleys, where we only have a few level lines to reconstruct the map from, cf. \cref{fig:mountain2}. Later we show how the reconstruction improves as the level lines increase, cf. \cref{fig:mountain_sequence}.
    
    \begin{figure}[!htbp]
        \captionsetup[subfigure]{
                                justification=centering,
                                labelformat=empty,
                                }
        \centering
        \begin{minipage}{1.0\textwidth}
            \centering
            \subfloat[][Ground truth]
                {
                \includegraphics[width=0.30\linewidth,
                                height=0.28\linewidth,
                                angle=0]{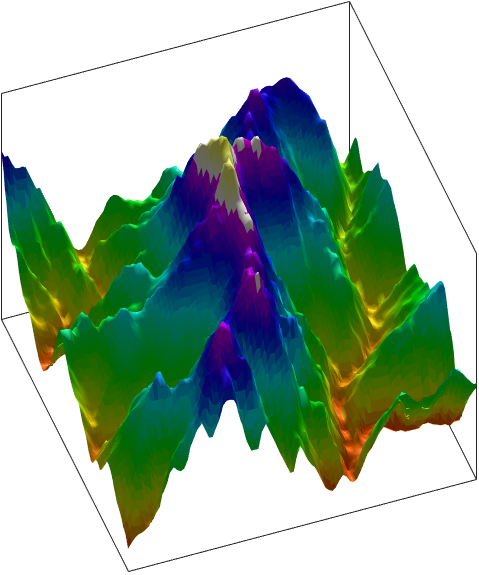}
                \label{subfig:mt_gt_1}
                }
            \hspace{0.01\linewidth}
            \subfloat[][$| \D^1 I |_F$ regularizer]
                {
                \includegraphics[width=0.30\linewidth,
                                height=0.28\linewidth,
                                angle=0]{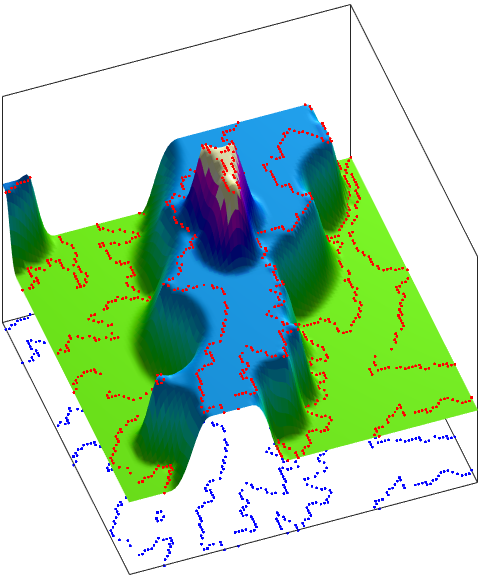}
                \label{subfig:mt_tv1}
                }
        \end{minipage}
        \begin{minipage}{1.0\textwidth}
            \centering
            \subfloat[][$| \D^2 I |_F$ regularizer]
                {
                \includegraphics[width=0.30\linewidth,
                                height=0.28\linewidth,
                                angle=0]{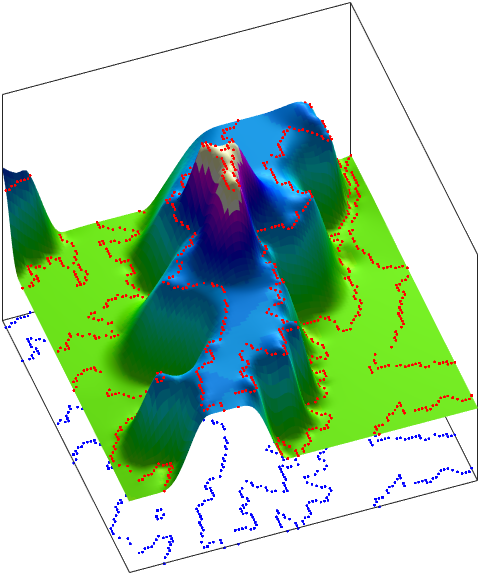}
                \label{subfig:mt_tv2}
                }
            \hspace{0.01\linewidth}
            \subfloat[][Anisotropic third order regularizer of \cite{LMS2013}]
                {
                \includegraphics[width=0.30\linewidth,
                                height=0.28\linewidth,
                                angle=0]{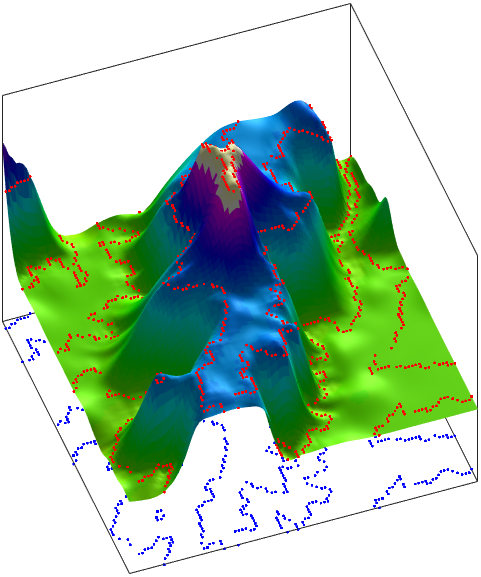}
                \label{subfig:mt_aniso_tv3}
                }
            \vspace{-0.01\linewidth}
        \end{minipage}
        \caption{Illustrating the challenges in reconstruction a real 3D map (mountain and valley) with only few level lines. The reconstructions here are from only three level lines (in blue on the plane, in red on the surface). The first two reconstructions, the upper-right and the lower-left sub-figure, are based on the first order and the second order total variation regularizer, respectively, but without any vector matching. The last reconstruction (the lower-right sub-figure) is based on the anisotropic third order regularizer proposed in \cite[$R^{(3)}_1(u)$ on page 5]{LMS2013}.}
        \label{fig:mountain1}
    \end{figure}

    To see the effectiveness of our model we first test it without the vector matching ($\a = 0$), cf. \cref{fig:mountain1}. In this test there are only three level lines available (blue lines on the floor or red lines on the surface), where the height values are given. The first two extensions are using the first order total variation $| \D^1 I |_F$ and the second order total variation $| \D^2 I |_F$, respectively. The reconstructed surfaces in both runs show difficult to reconstruct the valley in the ground truth. $\te=10^5$, $g=0$ and $h=10$ for the top right sub-figure, $g=10$ and $h=0$ for the bottom left sub-figure. The first one shows typical staircase effect while the second one favors slope. However in the valley lacking of information the only possible solution is the flattened structure. This is what we have already seen in previous experiments.

    {\sloppy In the last sub-figure of \cref{fig:mountain1}, we present the reconstruction using the anisotropic third order regularizer model of \cite[$R^{(3)}_1(u)$ in model (4) on page 5]{LMS2013}. The reconstructed surface improves a bit area still fails in reconstructing the valley, we believe it is because the level lines are too few to have enough similarity between the level lines which the algorithm requires, making it difficult to reconstruct the valley. For the experiment, we have solved the problem in \cite{LMS2013} using our algorithm based on the augmented Lagrangian.}

    \begin{figure}[!htbp]
        \captionsetup[subfigure]{
                                justification=centering,
                                labelformat=empty,
                                }
        \centering
        \begin{minipage}{1.0\textwidth}
            \centering
            \subfloat[][Ground truth]
                {
                \includegraphics[width=0.30\linewidth,
                                height=0.28\linewidth,
                                angle=0]{gt_3_1_x_plus}
                \label{subfig:mt_gt_2}
                }
            \hspace{0.01\linewidth}
            \subfloat[][The sign($\mbf{v}_{\G}$) is based on the isotropic surface]
                {
                \includegraphics[width=0.30\linewidth,
                                height=0.28\linewidth,
                                angle=0]{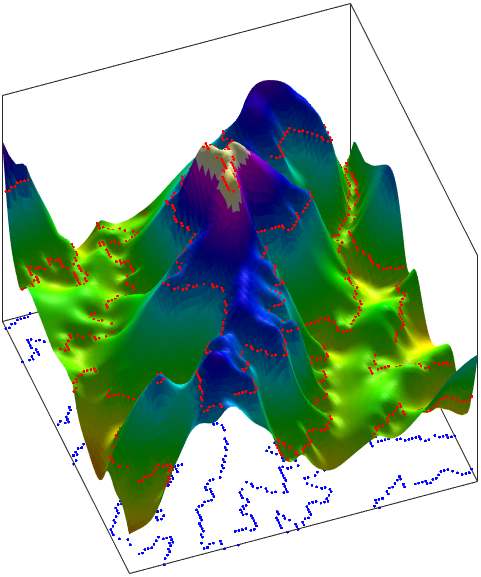}
                \label{subfig:mt_vct_w}
                }
        \end{minipage}
        \begin{minipage}{1.0\textwidth}
            \centering
            \subfloat[][The sign($\mbf{v}_{\G}$) is chosen adaptively]
                {
                \includegraphics[width=0.30\linewidth,
                                height=0.28\linewidth,
                                angle=0]{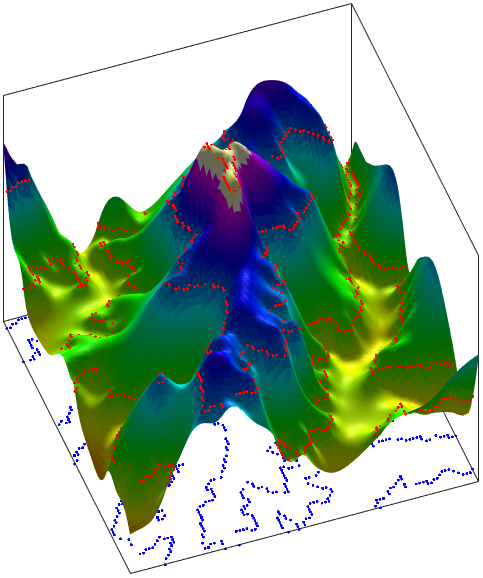}
                \label{subfig:mt_vct_r}
                }
            \hspace{0.01\linewidth}
            \subfloat[][The sign($\mbf{v}_{\G}$) is based on the ground truth]
                {
                \includegraphics[width=0.30\linewidth,
                                height=0.28\linewidth,
                                angle=0]{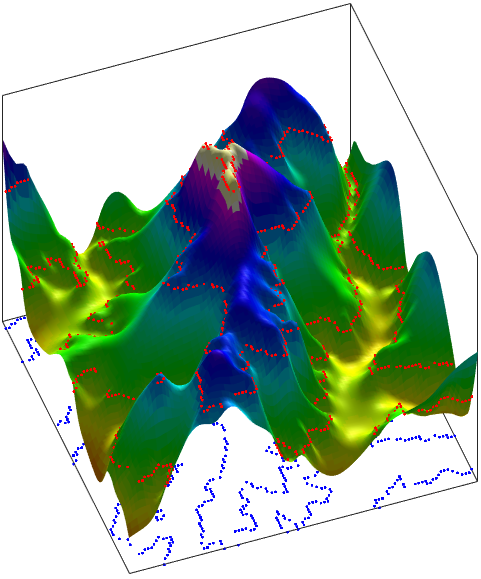}
                \label{subfig:mt_vct_gt}
                }
            \vspace{-0.01\linewidth}
        \end{minipage}
        \caption{Reconstructions of the same real 3D map (mountain and valley) as in \cref{fig:mountain1}, from its three level lines (blue lines on the floor or red lines on the surface), using our proposed model (\ref{eqs:opt_normal}). The different sub-figures here correspond to the different choices of sign for the unit normal vector $\mbf{v}_{\G}$ in the model.}
        \label{fig:mountain2}
    \end{figure}

    We now apply our model (\ref{eqs:opt_normal}) involving the normal vector matching term, on the three level lines. The results are shown in \cref{fig:mountain2}. The unit normal vector $\mbf{v}_{\G}$ needed in the model are extracted from the level line $C_{\G}$. Two possible choices  $\hat{\mbf{v}}_{\G} = \mbf{v}_{\G} \cdot | \hat{\mbf{v}}_{\G} |= \pm \partial_{t}^{\perp} C_{\G}$  where $\partial_{t} C_{\G}$, with $C_{\G}$ being the parametrization of the level line, denotes the tangent vector along $C_{\G}$ obtained via finite differences between discrete points along the level line.

    In the reconstruction, the correct sign of the unit normal vector $\mbf{v}_{\G}$ is not known in advance. To determine the sign we first solve the isotropic model, i.e. the model (\ref{eqs:opt_normal})) with $\a=0$, in order to get the $\D I$. We then look at the value of $\varrho := \mbf{v}_{\G} \cdot \D I$. If $\varrho<0$, we set $\mbf{v}_{\G} = - \mbf{v}_{\G}$. The surface presented in the upper-right corner in \cref{fig:mountain2} has been reconstructed in this way, showing already a significant improvement over the reconstructions presented in \cref{fig:mountain1}, where we see a more accurate representation of the mountain, as well as the two rivers in the image. 
    
    This is further improved, as we update the sign($\mbf{v}_{\G}$) gradually instead of everywhere at once, cf. the figure in the lower-left corner of \cref{fig:mountain2}. We let $\Gamma_{\varepsilon}$ be the subset of $\Gamma$ where we update the sign($\mbf{v}_{\G}$). In each iteration, we look at the absolute value of the directional derivative, i.e. $\varepsilon(x)=|\mbf{v}_{\G}(x) \cdot \D I(x)|$, and compare it with a threshold $\varepsilon_{threshold}$. Starting with the empty set $\G_{\varepsilon} = \varnothing$, we gradually include parts of $\Gamma$ whenever $\varepsilon(x)$ is found to be larger than the threshold in those parts. In other words, if $\varepsilon(x)>\varepsilon_{threshold}$ on $\G$ then update $\G_{\varepsilon} \leftarrow \G_{\varepsilon}\cup \{ x\in \G:\varepsilon(x) > \varepsilon_{threshold} \}$. At start we have the isotropic case, i.e. $\a=0$. $\a\neq0$ on $\G_{\varepsilon}$. If $\varrho(x)<0$ on $x\in\G_{\varepsilon}$ then we set $\mbf{v}_{\G} = - \mbf{v}_{\G}$ at that point. As the iteration continues, we expect $\G_{\varepsilon} = \G$. Gradually updating the level lines has given even better reconstruction which is very close to what we could achieve if sign($\mbf{v}_{\G}$) is extracted from the ground truth, cf. the figure in the lower-right corner of \cref{fig:mountain2}).

    As we can see, even with a few number of level lines, we have an almost perfect reconstruction. However, it becomes better as we increase the number of level lines, cf. \cref{fig:mountain_sequence}.

    In figures \ref{fig:mountain2}--\ref{fig:mountain_sequence}, the parameters for the regularizers has been $g=1$ and $h=0$, the data fidelity parameter has been $\te=10^5$, and the matching term parameter has been $\a=30$. In addition, $h=0.1$ has been on the boundary. $\varepsilon_{threshold} = 0.2$ has been used to determine the sign($\mbf{v}_{\G}$).

\subsubsection*{Real map II}
\begin{figure}[!htbp]
        \captionsetup[subfigure]{
                                justification=centering,
                                labelformat=empty,
                                }
        \centering
        \begin{minipage}{1.0\textwidth}
            \centering
            \subfloat[][Anisotropic third order model of \cite{LMS2013}]
                {
                \includegraphics[width=0.25\linewidth,
                                height=0.25\linewidth,
                                angle=0]{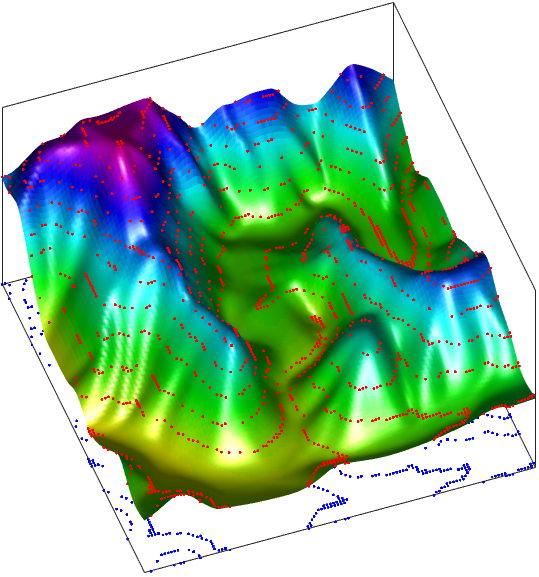}
                \label{subfig:cam_3_cam_n}
                }
            \hspace{0.01\linewidth}
            \subfloat[][Proposed model]
                {
                \includegraphics[width=0.25\linewidth,
                                height=0.25\linewidth,
                                angle=0]{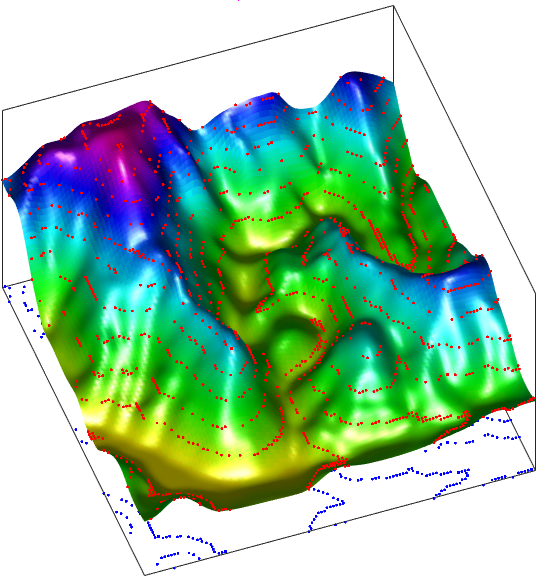}
                \label{subfig:cam_3_vct}
                }
            \hspace{0.01\linewidth}
            \subfloat[][Ground truth]
                {
                \includegraphics[width=0.25\linewidth,
                                height=0.25\linewidth,
                                angle=0]{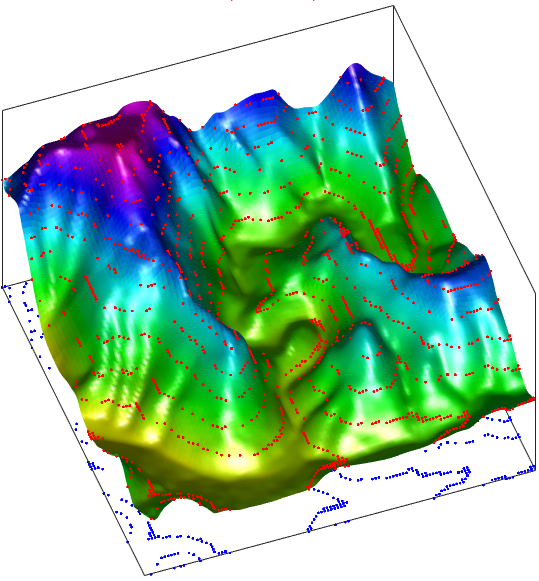}
                \label{subfig:cam_3_gt}
                }
        \end{minipage}
        \begin{minipage}{1.0\textwidth}
            \centering
            \subfloat[][Anisotropic third order model of \cite{LMS2013}]
                {
                \includegraphics[width=0.25\linewidth,
                                height=0.25\linewidth,
                                angle=0]{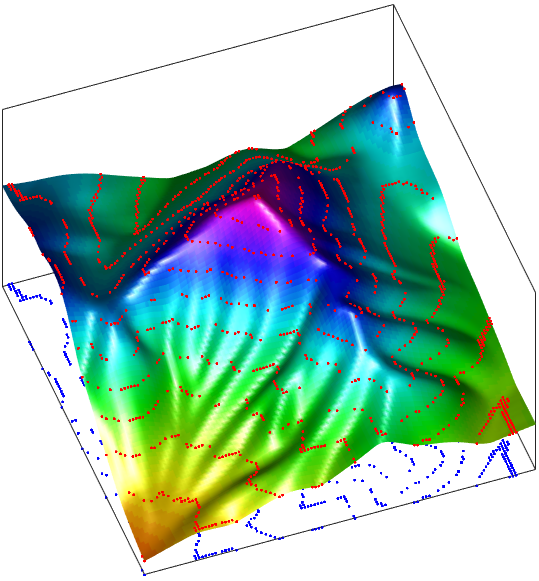}
                \label{subfig:cam_4_cam_n}
                }
            \hspace{0.01\linewidth}
            \subfloat[][Proposed model]
                {
                \includegraphics[width=0.25\linewidth,
                                height=0.25\linewidth,
                                angle=0]{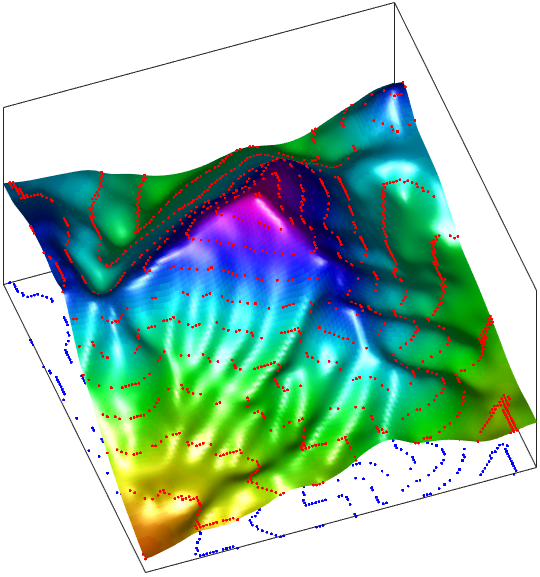}
                \label{subfig:cam_4_vct}
                }
            \hspace{0.01\linewidth}
            \subfloat[][Ground truth]
                {
                \includegraphics[width=0.25\linewidth,
                                height=0.25\linewidth,
                                angle=0]{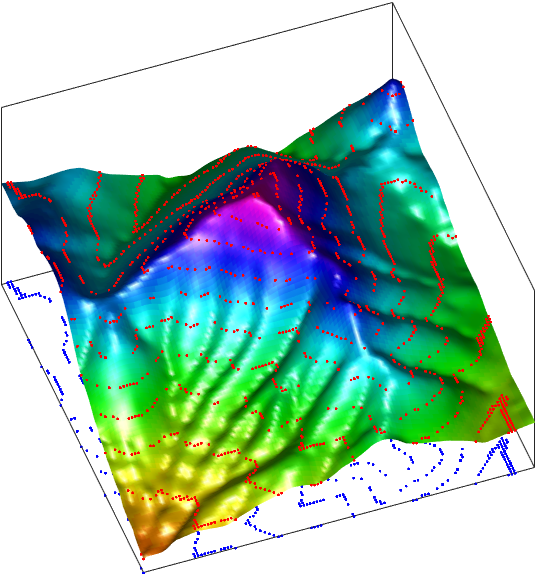}
                \label{subfig:cam_4_gt}
                }
            \vspace{-0.01\linewidth}
        \end{minipage}
        \caption{Comparing the proposed model with the anisotropic third order model proposed in \cite[$R^{(3)}_1(u)$ on page 5]{LMS2013} on two different real 3D maps from the same paper, one map in each row, with eleven level lines.}
        \label{fig:examples_from_cam}
    \end{figure}

    In the experiment we compare our model with the third order anisotropic model of \cite{LMS2013}. For the experiment we use two of the examples from \cite{LMS2013}, reproducing their results using our algorithm based on the augmented Lagrangian. The results are presented in \cref{fig:examples_from_cam}, one example in each row. For our model, $\te = 10^5$ (data fidelity), $\a = 10$ (vector matching) and $g = 0.1$ (the second order regularizer). For both, $h = 0$ inside and $h = 2$ on the boundary. $\varepsilon_{threshold} = 0.1$ has been used to determine sign($\mbf{v}_{\G}$). The parameters for augmented Lagrangian are set as $c_Q = 20$, $c_P = 1$ and $c_E = 50$ in this experiment. As seen from the figure, with enough level lines, both models perform well in capturing the anisotropy, however, our model managed to capture even the small variations, details like the small hill top etc.. This is because the model uses the geometry of the level lines.

\subsubsection*{Real map III}
\begin{figure}[!htbp]
        \captionsetup[subfigure]{
                                justification=centering,
                                labelformat=empty,
                                }
        \centering
        \begin{minipage}{1.0\textwidth}
            \centering
            \subfloat[][Ground truth]
                {
                \includegraphics[width=0.30\linewidth,
                                height=0.28\linewidth,
                                angle=0]{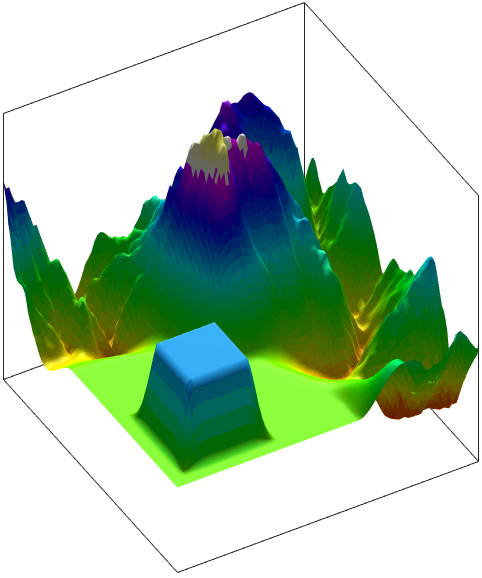}
                \label{subfig:mtb_gt}
                }
            \hspace{0.01\linewidth}
            \subfloat[][$| \D^2 I |_F$]
                {
                \includegraphics[width=0.30\linewidth,
                                height=0.28\linewidth,
                                angle=0]{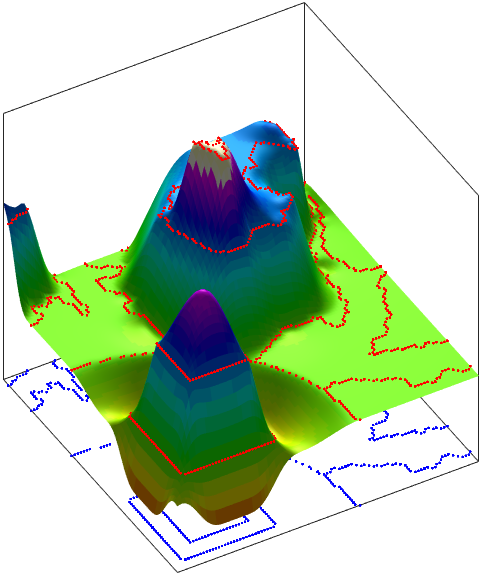}
                \label{subfig:mtb_tv2}
                }
        \end{minipage}
        \begin{minipage}{1.0\textwidth}
            \centering
            \subfloat[][$| \D^1 I |_F$]
                {
                \includegraphics[width=0.30\linewidth,
                                height=0.28\linewidth,
                                angle=0]{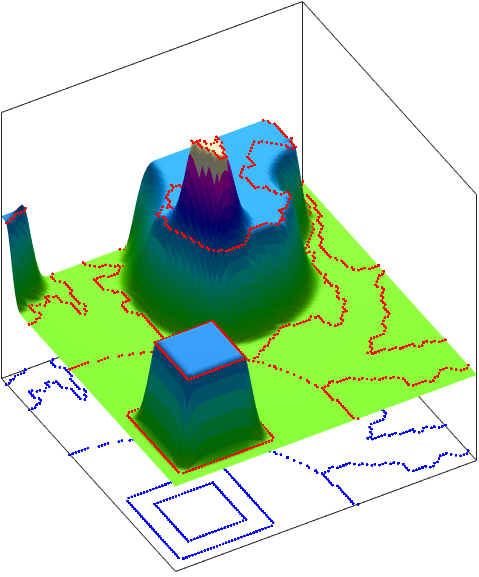}
                \label{subfig:mtb_tv1}
                }
            \hspace{0.01\linewidth}
            \subfloat[][Proposed model]
                {
                \includegraphics[width=0.30\linewidth,
                                height=0.28\linewidth,
                                angle=0]{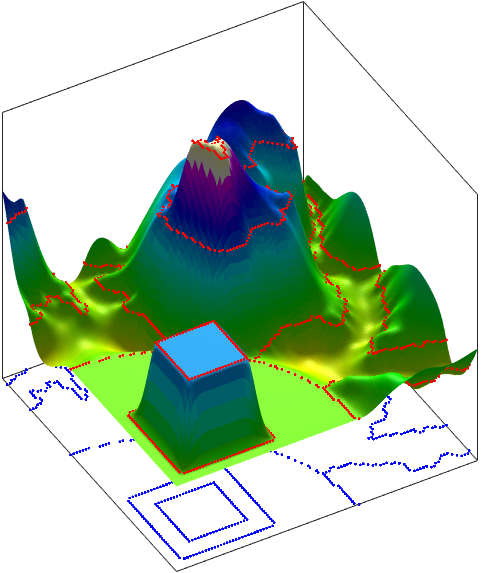}
                \label{subfig:mtb_vct}
                }
            \vspace{-0.01\linewidth}
        \end{minipage}
        \caption{
        Illustrating the effectiveness of the proposed model (\ref{eqs:opt_normal}) in reconstructing real 3D map with mixed structures, e.g. mountains, valleys, and buildings, from few level lines. The reconstructions here correspond to three level lines.}
        \label{fig:architecture}
    \end{figure}

    In this experiment we apply our model to reconstruct mixed structures with regular and irregular geometries. The top left sub-figure of \cref{fig:architecture} is the ground truth, which we reconstruct from its three level lines. Just as before, we apply our model with or without the vector matching and we observe the clear importance of having the vector matching. For the experiment, $\te = 10^5$ (data fidelity) on $\S$ and $h=0.1$ at boundary. For the top right sub-figure $g=10^{-2}$ and $h=0$ while $g=0$ and $h=50$ for the bottom left sub-figure. In the last sub-figure, $g=10^{-2}$ for the area with mountain, $h=50$ and $h=5$ for the flat base and tip of the building, respectively.  $\a=10$ and $\varepsilon_{threshold} = 0.2$ has been used to determine sign($\mbf{v}_{\G}$).

\section{Conclusions}
\label{sec:conclusions}

We have proposed a model consisting of a vector matching term (normal vector matching) to account for the anisotropy, together with a first-order and a second-order total variation regularizer (isotropic) under a fidelity constraint on the height. The model is able to effectively capture the irregularity along the given level lines, e.g., kinks or creases along level lines, and is able to recover surfaces from only a small number of level lines. 
    
We have proposed an effective way, based on the augmented Lagrangian method, to solve the model. In our algorithm, each sub-problem has either a closed form solution or a fast solver. We have derived closed form solution for the minimization problems containing $L^1$ term, provided a simple approach to derive the solution. For the inhomogenous modified Helmholtz equation (IMHE) with constant coefficient, we have given a closed form solution based on the discrete cosine transform. For the IMHE with variable coefficient (scalar function), we use the preconditioned conjugate gradient method, using the simplest yet effective preconditioner, the diagonal preconditioner.
    
The choice of parameters in our model plays a crucial role in the reconstruction. This choice may be turned automatically through some machine learning algorithm. This is a topic of future work.

\section*{Acknowledgments}
We thank Carola-Bibiane Sch\"onlieb and Jan Lellmann for the fruitful discussions in the beginning of this project, and particularly to Jan Lellmann for verifying our reconstruction in \cref{fig:mountain1} (the lower-right sub-figure) using their model.

\bibliographystyle{siamplain}
\bibliography{refs}

\end{document}